\documentclass{amsart}

\font\esc=cmcsc10 scaled 800

\let\what=\widehat
\let\wtil=\widetilde
\let\sse=\subseteq
\let\noi=\noindent
\let\vphi=\varphi
\let\veps=\varepsilon
\let\limply=\Longrightarrow

\def\span{{\kern.5pt{\rm span}\kern1pt}}
\def\smallfrac#1#2{{\textstyle{\frac{#1}{#2}}}}

\def\conv{{\;\longrightarrow\;}}
\def\wconv{{{\buildrel_{\scriptstyle w}\over\conv}}}
\def\sconv{{{\buildrel_{\scriptstyle s}\over\conv}}}

\def\query{{{\buildrel_{^{\scriptstyle ?}}\over\limply}}}

\def\jconv{{\conv_{_{\kern-13.5pt\scriptstyle j\kern9pt}}}}
\def\kconv{{\conv_{_{\kern-13.5pt\scriptstyle k\kern9pt}}}}

\def\emap{\hbox to25pt{\rightarrowfill}}

\def\nmap{\Big\uparrow}

\def\dsse{{\kern7pt{\overline{\phantom{\overline:}}{\kern-7pt\sse}}\kern3pt}}

\def\A{{\mathcal A}}
\def\B{{\mathcal B}}
\def\C{{\kern.5pt\mathcal C}}
\def\H{{\mathcal H}}
\def\M{{\mathcal M}}
\def\Oe{{\mathcal O}}
\def\X{{\mathcal X}}

\def\BX{{\B[\X]}}
\def\BH{{\B[\H]}}

\def\CC{{\mathbb C\kern.5pt}}
\def\FF{{\mathbb F\kern.5pt}}
\def\DD{{\mathbb D\kern.5pt}}
\def\NN{{\mathbb N\kern.5pt}}
\def\RR{{\mathbb R\kern.5pt}}
\def\TT{{\mathbb T\kern.5pt}}
\def\ZZ{{\mathbb Z\kern.5pt}}

\newsymbol\varnothing 203F
\let\void=\varnothing

\def\matrix#1{\null\,\vcenter{
		\normalbaselines\mathsurround=0pt\ialign{
		\hfil $##$
		\hfil && \quad
		\hfil $##$
		\hfil \crcr
		\mathstrut \crcr
		\noalign{\kern-\baselineskip}#1 \crcr
		\mathstrut \crcr
		\noalign{\kern-\baselineskip} \crcr }}\,}

\def\newmatrix#1{\null\,\vcenter{
		\baselineskip=8pt\mathsurround=-0pt\ialign{
		\hfil ${##}$
		\hfil &&
		\hfil ${##}$
		\hfil \crcr
		\mathstrut \crcr
		\noalign{\kern-\baselineskip}#1 \crcr
		\mathstrut \crcr
		\noalign{\kern-\baselineskip} \crcr }}\!}

\begin{document}

\vglue-64pt\noindent
\hfill{\it Results in Mathematics}\/,
{\bf 79}(5) (2024) 1--30

\vglue20pt
\title[Weak Supercyclicity]
      {Weak Supercyclicity --- An Expository Survey}
\author[Carlos Kubrusly]
       {Carlos S. Kubrusly}
\address{Catholic University of Rio de Janeiro, Brazil}
\email{carlos@ele.puc-rio.br}
\renewcommand{\keywordsname}{Keywords}
\dedicatory{Dedicated to the memory of my PhD supervisor Ruth F$.$ Curtain
            $($1941--2018$\kern1pt)$}
\date{October 12, 2023; revised: March 25, 2024}

\begin{abstract}
This is an exposition on supercyclicity and weak supercyclicity, especially
designed to advance further developments in weakly supercyclicity, which is a
recent research field showing significant momentum during the past two
decades$.$ For operators on a normed space, the present paper explores the
relationship among supercyclicity and three versions of weak supercyclicity,
namely, weak l-sequential, weak sequential, and weak supercyclicities, and
their connection with strong stability, weak stability, and weak
quasistability$.$ A survey of the literature on weak supercyclicity is
followed by an analysis of the interplay between supercyclicity and strong
stability, as well as between weak supercyclicity and weak stability$.$ A
description of the spectrum of weakly l-sequentially supercyclic operators is
also given.
\end{abstract}

\maketitle

\vskip-15pt\noi
\section{Introduction}

This is an expository-survey paper focusing on weak forms of supercyclicity
and their connection with weak stability of bounded linear operators on a
normed space$.$ It has been devised to offer an organised, summarised, and
unified collection of results in such a recent research area, as well as to
exhibit some motivating open questions on it$.$ The target is weak
supercyclicity, in particular, weak l-sequential super\-cyclicity, whose
relatively recent mathematical developments began around 2005, and have
become a particularly active research field since then.

\vskip6pt
It has been known for some time that supercyclicity in the norm topology, also
referred to as strong supercyclicity, implies strong stability$.$ This means
that if the projective orbit of a power bounded linear operator $T$ on a
normed space $\X$ is dense in $\X$ for some vector $y$ in $\X$, then
${T^nx\to0}$ for every ${x\in\X}$, with both density and convergence in the
norm topology$.$ Symbolically: ${\Oe_T([y])^-\!=\X\limply T^n\!\sconv O}.$
This is an important result from \cite[Theorem 2.2]{AB}, which becomes a major
motivation behind the quest for a possible relationship between some form of
weak supercyclicity and weak stability$.$ Does such an implication survive the
transition from norm topology to weak topology$?$ In fact, it was asked in
\cite{KD1} whether a stronger version of weak supercyclicity, viz., weak
l-sequential supercyclicity, implies weak stability$.$ As discussed in
Section 3, weak l-sequential supercyclicity is weaker than supercyclicity in
the norm topology but stronger than supercyclicity in the weak sequential
topology, and so stronger than supercyclicity in the weak topology.

\vskip6pt
The present paper brings an exposition on weak supercyclicity in general,
including a survey of the recent literature on weak supercyclicity and, in
particular, a discussion on the above mentioned question.

\vskip6pt
All terms used above will be defined here in due course$.$ The paper is
organised into six more sections$.$ Since it deals with many aspects related
to linear dynamics, involving a large amount of specific jargon and
concepts, a section on notation and terminology will be split into two parts
in order to make it easier to read$.$ Basic no\-tation and terminology are
posed in Section 2, while notions involving cyclicity, hy\-percyclicity,
supercyclicity, and three forms of weak supercyclicity are considered in
Section 3$.$ That being said, the literature specifically on forms of
{\it weak supercyclicity}\/ is surveyed in Section 4$.$ Such literature is
classified according to the main target of each paper, roughly in
chronological order, \hbox{into eight subgroups, namely},
\vskip4pt
(1) pre-weak supercyclicity,
\vskip3pt
(2) the beginning of weak supercyclicity,
\vskip3pt
(3) further forms of weak supercyclicity,
\vskip3pt
(4) clarifying three forms of weak supercyclicity,
\vskip3pt
(5) weak supercyclicity for operators close to normal,
\vskip3pt
(6) weak supercyclicity for compact operators,
\vskip3pt
(7) weak supercyclicity for composition operators,
\vskip3pt
(8) weak supercyclicity and weak stability.
\vskip5pt\noi
The fundamental arguments ensuring that strong supercyclicity implies strong
stability are considered in detail in Section 5$.$ The special themes of weak
supercyclicity (in particular, weak l-sequential supercyclicity) and weak
stability are reviewed in Chapter 6$.$ We close the paper by proving a
characterisation of the spectrum of weakly l-sequentially supercyclic
operators in Section 7.
\vskip-4pt\noi

\section*{{\bf 2}. Notation and Terminology, Part I: Basic Notions}

Throughout the paper, linear spaces are over a scalar field $\FF\kern-1pt$,
which here is always $\RR$ or $\CC.$ Let $\BX$ denote the normed algebra of
all bounded linear \hbox{transformations} of a normed space $\X$ onto
itself$.$ Elements of $\BX$ will be referred to as operators$.$ Both norms on
a normed space $\X$ and the induced uniform norm on $\BX$ will be denoted by
the same symbol ${\|\cdot\|}.$ Let $\X^*\!$ be the dual of a normed space
\hbox{$\X$, and let} ${T^*\!\in\B[\X^*]}$ stand for the normed-space adjoint
of ${T\kern-1pt\in\BX}.$ We use the same notation ${T^*\!\in\BH}$ for the
Hilbert-space adjoint if ${T\kern-1pt\in\BH}$ acts on a Hilbert space $\H.$
The null operator and the identity operator will be denoted by $O$ and $I$,
\hbox{respectively}.

\vskip6pt
We begin with a standard definition (see, e.g.,
\cite[Definition 1.13.2]{Meg}$\kern.5pt$).

\vskip6pt\noi
{\bf Definition.}$\kern-.5pt$
Let $\X$ be a normed space$.$ An $\X$-valued sequence $\{x_n\}$ converges
weak\-ly, or $\{x_n\}$ is weakly convergent, if the scalar-valued sequence
$\{f(x_n)\}$ converges for every ${f\in\X^*}\!$ in the sense that there exists
${z\in\X}$ for which the scalar-valued sequence $\{f(x_n)\}$ converges to
$f(z)$ for every ${f\in\X^*}\!.$

\vskip6pt
Therefore, an $\X$-valued sequence $\{x_n\}$ converges weakly if
${\lim_nf(x_n)\to f(z)}$ for some ${z\in\X}$, for every ${f\in\X^*}.$ Usual
notation: ${x_n\wconv z}$ or $z={w\hbox{\,-}\lim_n x_n}.$ An $\X$-valued
sequence $\{x_n\}$ is strongly convergent, or converges in the norm topology,
if ${\lim_n\|x_n-z\|\to0}$ for some ${z\in\X}.$ Usual notation$:$
${x_n\to z}$.

\vskip6pt
For any ${T\kern-1pt\in\BX}$, take its $\BX$-valued power sequence $\{T^n\}$
(with $n$ running over the nonnegative integers) and, for each ${x\in\X}$,
consider the $\X$-valued sequence $\{T^nx\}.$ An operator $T$ is power
bounded if ${\sup_n\|T^n\|<\infty}.$ If $\X$ is a Banach space, then
this is equivalent to ${\sup_n\|T^nx\|<\infty}$ for every ${x\in\X}$ by the
Banach Steinhaus Theorem$.$ Moreover, an operator $T$ is of class
$C_{1{\textstyle\cdot}}$ if
$$
\hbox{${T^nx\not\to0}\;$ for every nonzero ${x\in\X}$}
$$
\goodbreak\noi
(i.e., ${\|T^nx\|\not\to0}$ for every ${0\ne x\in\X}).$ The opposite notion is
that of strongly stable operators$.$ An operator $T$ is strongly stable, or
$T$ is of class $C_{0{\kern.5pt\textstyle\cdot}}\kern-1pt$, if
$$
\hbox{${T^nx\to0}\;$ for every ${x\in\X}$}
$$
(i.e., ${\|T^nx\|\to0}$ for every ${x\in\X}).$ On the other hand, $T$ is
weakly stable if
$$
\hbox{ ${T^nx\wconv0}\;$ for every ${x\in\X}$}.
$$
Equivalently, an operator $T$ is weakly stable if
$$
{\lim}_n|f(T^nx)|=0
\;\;\hbox{for every}\;
x\in\X,
\;\;\hbox{for every}\;
f\in\X^*\!.
$$
Usual notation for strong and weak stability$:$ ${T^n\!\sconv O}$ and
${T^n\!\wconv O}$, respectively.

\vskip6pt
In general, strong convergence implies weak convergence, and so, in
particular, strong stability implies weak stability$.$ If $\X$ is a Banach
space, then weak stability implies power boundedness by the Banach--Steinhaus
Theorem$.$ If $\X$ is a complex Banach space, then let $\sigma(T)$ denote the
spectrum of $T$ and let $r(T)=\sup_{\lambda\in\sigma(T)}|\lambda|$ stand for
the spectral radius of $T\!.$ Recall that ${0\le r(T)\le\\|T\|}$, and $T$ is
quasinilpotent or normaloid if ${r(T)=0}$ or ${r(T)=\|T\|}$, respectively$.$
The spectrum of a power bounded operator is included in the closed unit
disk$.$ Summing up:
$$
T^n\!\sconv O
\quad\limply\quad
T^n\!\wconv O
\quad\limply\quad
{\sup}_n\|T^n\|<\infty
\quad\limply\quad
r(T)\le1.
$$
\vskip0pt
A final note on basic notation: we write $\ell^p$ for $\ell^p(\ZZ)$ and
$\ell_+^p$ for $\ell^p(\NN)$ or $\ell^p(\NN_0)$, for any ${p\ge\kern-1pt1}$,
where $\ZZ$, $\NN$, and $\NN_0$ stand for the sets of all integers, all
positive integers, and all nonnegative integers, respectively.

\section*{{\bf 3}. Notation and Terminology, Part II: Notions of Cyclicity}

This section pulls together the specific notation and terminology required
to describe all kinds of cyclicity$.$ Before introducing the standard notions
of hypercyclicity, supercyclicity, and cyclicity, and also the variants of
weak supercyclicity, we need the notions of weak l-sequential denseness
and projective orbit.

\vskip6pt\noi
{\bf \small 3.1.}
{\sc weak l-sequential denseness}$.$
Let $A$ be a subset of a normed space $\X.$ Its closure in the norm topology
on $\X$ is denoted by $A^-\!$, and its weak closure in the weak topology on
$\X$ is denoted by $A^{-w}\!.$ Thus $A$ is dense (in the norm topology) or
weakly dense (in the weak topology) if $A^-\!=\X$ or $A^{-w}\!=\X$,
respectively$.$ Since the weak topology is not metrizable, the notions of weak
closure and weak sequential closure do not necessarily coincide$.$ A set $A$
is weakly sequentially closed if every $A$-valued weakly convergent sequence
has its limit in $A.$ Let the weak sequential closure of $A$ be denoted by
$A^{-ws}\!$, which is the smallest weakly sequentially closed subset of $\X$
including $A.$ That is, $A^{-ws}\!$ is the intersection of all weakly
sequentially closed subsets of $\X$ that include $A.$ So a set $A$ is weakly
sequentially dense if $A^{-ws}\!=\X.$ Differently, the {\it weak limit set}\/
$A^{-wl}$ of $A$ is the set of all weak limits of weakly convergent $A$-valued
sequences, that is,
$$
A^{-wl}=\big\{x\in\X\!:x=w\hbox{\,-}{\lim}_nx_n\;\hbox{with}\;\,x_n\in A\big\}.
$$
We say that $A$ is {\it weakly l-sequentially dense}\/ if $A^{-wl}\!=\X.$ It
is known that (see, e.g., \cite[Proposition 2.1]{Kub2}$\kern.5pt$)
\vskip4pt\noi
$$
A^-\!\sse A^{-wl}\!\sse A^{-ws}\!\sse A^{-w}\!.
$$
\vskip3pt\noi
These inclusions may be proper in general (see, e.g., \cite[pp.38,39]{Shk1},
\cite[pp.259,260]{BM2}, \cite[pp.10,11]{Kub2}$\kern.5pt).$ Consequently, if
$A$ is dense in the norm topology (i.e., ${A^-\!=\X}$), then it is dense with
respect to all notions of denseness defined above$.$ Recall that if a set $A$
is convex, then ${A^-\!=A^{-w}}$ (see, e.g.,
\cite[Theorem 2.5.16]{Meg}$\kern.5pt).$ Thus, if $A$ is convex, then the above
inclusions become identities, and so for a convex set all the above notions of
denseness coincide.

\vskip6pt
A vector ${x\in\X}$ is a {\it weak limit point}\/ of a set ${A\sse\X}$ if
${x\in A^{-wl}}\!$, and it is a {\it weak accumu\-lation point}\/ of $A$ if
${x\in A^{-w}}\!.$ Thus the above inclusions show that weak limit points
are weak accumulation points, but the converse, in general, fails.

\vskip6pt\noi
{\bf \small 3.2.}
{\sc projective orbit}$.$
Given an operator ${T\kern-1pt\in\kern-1pt\BX}$ on a normed space
$\X\kern-1pt$, consider again its power sequence $\{T^n\}$ (with $n$ running
over the nonnegative integers)$.$ The {\it orbit}\/ $\Oe_T(y)$ or
${\rm Orb}\kern1pt(T\kern-1pt,y)$ of a vector ${y\in\X}$ under $T$
is the set
$$
\Oe_T(y)={\bigcup}_nT^n\{y\}={\bigcup}_n\{T^ny\}
=\big\{T^ny\in\X\!:\hbox{for every nonnegative integer $n$}\big\}.
$$
The orbit $\Oe_T(A)$ of a set ${A\sse\X}$ under $T$ is likewise defined$:$
$\Oe_T(A)=\bigcup_nT^n(A)\!=\!\bigcup_{z\in A}\Oe_T(z).$ Let $[x]=\span\{x\}$
stand for the sub\-space of $\X$ spanned by a singleton $\{x\}$ at a vector
${x\in\X}$, which is a one-dimensional subspace of $\X$ whenever $x$ is
nonzero$.$ The {\it projective orbit}\/ of a vector ${y\in\X}$ under an
operator ${T\kern-1pt\in\kern-1pt\BX}$ is the orbit of span of $\{y\}$; that
is, the orbit $\Oe_T([y])$ of $[y]$:
$$
\Oe_T([y])
=\Oe_T(\span\{y\})
={\bigcup}_{z\in[y]}\Oe_T(z)
=\big\{\alpha T^ny\in\X\!:\,\alpha\in\FF,\;n\in\NN_0\big\}.
$$

\vskip2pt\noi
{\bf \small 3.3.}
{\sc hypercyclicity, supercyclicity, and cyclicity}$.$
It is natural and usual to classify the concepts of cyclicity into three
classes, namely, hypercyclicity, supercyclicity, and cyclicity$.$ According to
the type of density (cf$.$ Subsection 3.1) upon which these notions are
defined, a further split goes as follows.

\begin{description}
\item$\kern-5pt$(A)
A vector $y$ in $\X$ is a {\it hypercyclic vector}\/ for $T$ if its orbit
$\Oe_T(y)$ is dense in $\X$ (i.e., dense in the norm topology):
$$
\Oe_T(y)^-\!=\X.
$$
As the norm topology is metrizable, a nonzero vector $y$ in the normed space
$\X$ is a hypercyclic vector for $T$ if for every $x$ in $\X$ there is a
sequence with ele\-ments in $\Oe_T(y)$ converging to it$.$ That is, ${y\in\X}$
is a hypercyclic vector for $T$ if for every $x$ in $\X$ there is a sequence
$\{T^{n_k}\}_{k\ge0}$ (which depends on $x$) with entries from
$\{T^n\}_{n\ge0}$ such that
$$
T^{n_k}y\to x
\qquad(\hbox{i.e.,}\;\; \|T^{n_k}y-x\|\to0).
$$
An operator $T$ is a {\it hypercyclic operator}\/ if it has a hypercyclic
vector$.$ Clearly, a hypercyclic operator is not power bounded.

\vskip6pt
\item$\kern-5pt$(B)
A vector $y$ in $\X$ is a {\it supercyclic vector}\/ for $T$ if its projective
orbit $\Oe_T([y])$ is dense in $\X.$ That is, if the orbit of the span of
$\{y\}$ is dense in $\X$ (i.e., dense in the norm topology):
$$
\Oe_T([y])^-\!=\X
\qquad(\hbox{i.e.,}\quad
\Oe_T(\span\{y\})^-\!=\X).
$$
Again, as the norm topology is metrizable, a nonzero vector $y$ in $\X$ is a
super\-cyclic vector for an operator $T$ if for every $x$ in $\X$ there is an
$\FF$-valued sequence $\{\alpha_k\}_{k\ge0}$ such that, for some sequence
$\{T^{n_k}\}_{k\ge0}$ with entries from $\{T^n\}_{n\ge0}$, the sequence
$\{\alpha_kT^{n_k}y\}_{k\ge0}$ converges to $x$ (in the norm topology):
$$
\alpha_kT^{n_k}y\to x
\qquad(\hbox{i.e.,}\quad
\|\alpha_kT^{n_k}y-x\|\to0).
$$
The scalar sequence $\{\alpha_k\}_{k\ge0}$ consists of nonzero numbers if
${x\ne0}.$ Both sequences $\{\alpha_k\}_{k\ge0}$ and $\{T^{n_k}\}_{k\ge0}$
depend on $x$ for each $y.$ An operator $T$ is a {\it supercyclic operator}\/
if it has a supercyclic vector.
\end{description}
\vskip-2pt

\vskip6pt
The next two definitions are weak counterparts of the above convergence
criteria$.$

\vskip6pt
\begin{description}
\item$\kern-5pt$(C)
A nonzero vector $y$ in $\X$ is a {\it weakly l-sequentially hypercyclic
vector}\/ for $T$ if for every $x$ in $\X$ there exists a sequence
$\{T^{n_k}\}_{k\ge0}$ with entries from $\{T^n\}_{n\ge0}$ such that the
$\X$-valued sequence $\{T^{n_k}y\}_{k\ge0}$ converges weakly to $x$,
$$
T^{n_k}y\wconv x\qquad(\hbox{i.e.,}\quad
f(T^{n_k}y-x)\to0
\quad\hbox{for every}\quad
f\in\X^*).
$$
The sequence $\{T^{n_k}\}_{k\ge0}$ depends on $x$ for each $y$ but does not
depend on $f.$ This means the orbit $\Oe_T(y)$ of $y$ under $T\kern-1pt$ is
weakly l-sequentially dense in $\X$ in the sense that the weak limit set of
$\Oe_T(y)$ coincides with $\X$:
$$
\Oe_T(y)^{-wl}=\X.
$$
An operator $T$ is said to be a {\it weakly l-sequentially hypercyclic
operator}\/ if it has a weakly l-sequentially hypercyclic vector.

\vskip6pt
\item$\kern-5pt$(D)
A nonzero vector $y$ in $\X$ is a {\it weakly l-sequentially supercyclic
vector}\/ for $T$ if for every $x$ in $\X$ there exists an $\FF$-valued
sequence $\{\alpha_k\}_{k\ge0}$ such that, for some sequence
$\{T^{n_k}\}_{k\ge0}$ with entries from $\{T^n\}_{n\ge0}$, the $\X$-valued
sequence $\{\alpha_kT^{n_k}y\}_{k\ge0}$ converges weakly to $x.$ In other
words, if
$$
\alpha_kT^{n_k}y\wconv x
\qquad(\hbox{i.e.,}\quad
f(\alpha_kT^{n_k}y-x)\to0
\quad\hbox{for every}\quad
f\in\X^*).
$$
The scalar sequence $\{\alpha_k\}_{k\ge0}$ consists of nonzero numbers if
${x\ne0}$, and both sequences $\{\alpha_k\}_{k\ge0}$ and $\{T^{n_k}\}_{k\ge0}$
depend on $x$ for each $y$ but do not depend on $f.$ This means the projective
orbit $\Oe_T([y])$ of $y$ under $T\kern-1pt$ is weakly l-sequentially dense in
$\X$ in the sense that the weak limit set of $\Oe_T([y])$ coincides with $\X$:
$$
\Oe_T([y])^{-wl}=\X.
$$
An operator $T$ is a {\it weakly l-sequentially supercyclic operator}\/ if it
has a weakly l-sequentially supercyclic vector.
\end{description}

\vskip6pt
Before proceeding, it is convenient to comment on two important points$.$
Remark 1 below shows that weak l-sequential supercyclicity and weak
l-sequential hyper\-cyclicity {\it are not topological notions}\/$.$ In the
subsequent Remark 2, it is recalled that any form of cyclicity makes
nontrivial sense only when approaching points outside the projective orbit,
or outside the orbit, respectively.

\vskip6pt\noi
{\bf Remark 1.}$\kern-.5pt$ {\sc weak l-sequential supercyclicity}\/$.$
It was observed in \cite{Kub2} that, for any subset $A$ of $\X$, the map
${A\mapsto A^{-wl}}$ is not a topological \hbox{closure operation}$.$ Thus
weak l-sequential supercyclicity (and weak l-sequential hypercyclicity) are
not topological notions$.$ Indeed, as defined in (D) (and in (C)) above, weak
l-sequential supercyclicity (and weak l-sequential hypercyclicity) come as
weak versions of a norm convergence criterion, but not as denseness of the
projective orbit (or of the orbit) in any topology$.$ In particular, weak
l-sequential supercyclicity (and weak l-sequential hypercyclicity) do not mean
denseness of the projective orbit (or of the orbit) in the weak topology or in
the weak sequential topology$.$ Therefore, weak l-sequential supercyclicity
(and weak l-sequential hypercyclicity) are not weak topology notions$.$
Consequently, it is not the case of applying weak topology techniques when
dealing with weak \hbox{l-sequential supercyclicity} (or with weak
l-sequential hypercyclicity)$.$ For a discussion along this line, see also
\cite{Shk1}.

\vskip6pt\noi
{\bf Remark 2.}$\kern-.5pt$ {\sc sequence with entries from a sequence}\/$.$
Take a nonzero ${x\in\kern-1pt\X}\kern-1pt.$ If ${x\in\Oe_T([y])}$ (or
${x\in\Oe_T(y)}$\/) for some vector $y$, then the sequence
$\{T^{n_k\!}\}_{k\ge0}$ mentioned in (A) to (D) can be thought of as a
constant infinite sequence or, equivalently, as a one-entry, thus finite,
subsequence of $\{T^n\}_{n\ge0}.$ For instance, in the supercyclicity cases
described in (B) and (D), if ${x\in\Oe_T([y])}$, then ${x=\alpha_0 T^{n_0}y}$
for some ${\alpha_0\ne0}$ and some ${n_0\ge0}.$ In such cases, both notions of
${\alpha_kT^{n_k}y\wconv x}$ and ${\alpha_kT^{n_k}y\to x}$, associated with
weak l-sequential supercyclicity and supercyclicity, coincide, where
convergence of the constant or finite sequence $\{\alpha_kT^{n_k}y\}$ means
that a limit is eventually reached$.$ Thus, in general, we use the expression
$$
\hbox{\it $\{T^{n_k}\}_{k\ge0}$ is a sequence with entries from\/
$\{T^n\}_{n\ge0}$}.
$$
However, if ${x\not\in\Oe_T([y])}$ (or ${x\not\in\Oe_T(y)}\kern.5pt$), this
means
$$
\hbox{\it $\{T^{n_k}\}_{k\ge0}$ is a subsequence of\/ $\{T^n\}_{n\ge0}$}.
$$
\vskip-2pt

\vskip6pt
Since the weak topology is not metrizable, the notions of weak closure and
weak sequential closure do not necessarily coincide.
\vskip4pt
\begin{description}
\item$\kern-5pt$(E)
A vector $y$ in $\X$ is a {\it weakly sequentially hypercyclic vector}\/ for
$T$ if its orbit $\Oe_T(y)$ is weakly sequentially dense in $\X$ (i.e., dense
in the weak sequential topology):
\vskip-4pt\noi
$$
\Oe_T(y)^{-ws}\!=\X.
$$
An operator $T$ is a {\it weakly sequentially hypercyclic operator}\/ if it
has a weakly sequentially hypercyclic vector$.$

\vskip4pt
\item$\kern-5pt$(F)
A vector $y$ in $\X$ is a {\it weakly sequentially supercyclic vector}\/ for
$T$ if its projective orbit $\Oe_T([y])$ is weakly sequentially dense in $\X$
(i.e., dense in the weak sequential topology):
\vskip-2pt\noi
$$
\Oe_T([y])^{-ws}\!=\X.
$$
An operator $T$ is a {\it weakly sequentially supercyclic operator}\/ if it
has a weakly sequentially supercyclic vector$.$

\vskip4pt
\item$\kern-5pt$(G)
A vector $y$ in $\X$ is a {\it weakly hypercyclic vector}\/ for $T$ if its
orbit $\Oe_T(y)$ is weakly dense in $\X$ (i.e., dense in the weak topology):
$$
\Oe_T(y)^{-w}\!=\X.
$$
An operator $T$ is a {\it weakly hypercyclic operator}\/ if it has a weakly
hypercyclic vector$.$
\vskip-2pt

\vskip4pt
\item$\kern-5pt$(H)
A vector $y$ in $\X$ is a {\it weakly supercyclic vector}\/ for $T$ if its
projective orbit $\Oe_T([y])$ is weakly dense in $\X$ (i.e., dense in the weak
topology):
$$
\Oe_T([y])^{-w}\!=\X.
$$
An operator $T$ is a {\it weakly supercyclic operator}\/ if it has a weakly
supercyclic vector$.$
\vskip-2pt

\vskip4pt
\item$\kern-5pt$(I)
A vector $y$ in $\X$ is a {\it cyclic vector}\/ for an operator $T$ if the
span of its orbit $\span\,\Oe_T(y)$ is dense in $\X$ (i.e., dense in the norm
topology):
$$
\big(\span\,\Oe_T(y)\big)^-\!=\X.
$$
An operator $T$ is a {\it cyclic operator}\/ if it has a cyclic vector.
\end{description}

\vskip6pt\noi
{\bf \small 3.4.}
{\sc The relation among hypercyclicity, supercyclicity, and cyclicity}$.$
Since $\span\,\Oe_T(y)$ is a convex subset of a normed space, its (norm)
closure coincides with its weak closure,
${(\span\,\Oe_T(y))^-\!=(\span\,\Oe_T(y))^{-w}}\!$, and hence weak cyclicity
boils down to cyclicity, and so do the intermediate concepts:
$$
\big(\span\,\Oe_T(y)\big)^-\!
=\big(\span\,\Oe_T(y)\big)^{-wl}\!
=\big(\span\,\Oe_T(y)\big)^{-ws}\!
=\big(\span\,\Oe_T(y)\big)^{-w}.
$$
\vskip-2pt

\vskip6pt
The above notions are related as follows.
\vskip6pt\noi
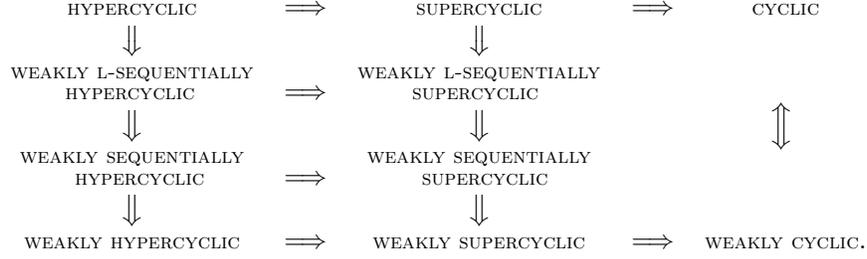
$$
\newmatrix{
& \hbox{\esc hypercyclic}
& \limply
& \hbox{\esc supercyclic}
& \limply
& \hbox{\esc cyclic}                                                   \cr
& \big\Downarrow
&
& \big\Downarrow
&
&                                                       \phantom{\Big|} \cr
& \vbox{\hbox{\esc weakly l-sequentially}
       \hbox{\esc \phantom{------} hypercyclic}}
& \limply
& \vbox{\hbox{\esc weakly l-sequentially}
       \hbox{\esc \phantom{------} supercyclic}}
&
&                                                                       \cr
& \big\Downarrow
&
& \big\Downarrow
&
& \Big\Updownarrow                                      \phantom{\Big|} \cr
& \vbox{\hbox{\esc weakly sequentially}
        \hbox{\esc \phantom{------} hypercyclic}}
& \limply
& \vbox{\hbox{\esc weakly sequentially}
        \hbox{\esc \phantom{------} supercyclic}}
&
&                                                                       \cr
& \big\Downarrow
&
& \big\Downarrow
&
&                                                       \phantom{\Big|} \cr
& \hbox{\esc weakly hypercyclic}
& \;\;\;\limply\;\;\;
& \hbox{\esc weakly supercyclic}
& \;\;\;\limply\;\;\;
& \hbox{\esc weakly cyclic}.                                            \cr}
$$

\vskip6pt
It is known that there is no hypercyclic operator on a finite-dimensional
normed space (see, e.g.,\cite[Proposition 1.1]{BM2} and
\cite[Theorem 2.58]{GP}$\kern.5pt$), although every nonzero vector is
trivially supercyclic for every operator on a one-dimensional space$.$
\hbox{Recall} that weak and norm topologies coincide in a finite-dimensional
space, and also that cyclicity implies separability for any normed space, as
it is spanned by a countable orbit, and therefore all forms of cyclicity in
the above diagram imply separability$.\kern-1pt$ This paper focuses mainly on
three forms of weak supercyclicity, and so all normed spaces here are
supposed to be infinite-dimensional and separable.

\vskip6pt
For basic properties concerning these notions of cyclicity, see, e.g.,
\cite[pp.38,39]{Shk1}, \cite[pp.259,260]{BM2}, \cite[pp.159,232]{GP},
\cite[pp.50,51]{KD1}$.$ In particular, questions related to some converses
in the above diagram asking, for instance, whether weak l-sequential
supercyclicity differs from weak sequential supercyclicity,
$$
\hbox{\it is there a weakly sequentially supercyclic not weakly
l-sequentially supercyclic}\,?
$$
as well as whether weak l-sequential hypercyclicity differs from weak
sequential hyper\-cyclicity, have been addressed in \cite[pp.71,72]{Shk1}$.$
However, by omitting the four entries in the above diagram containing the term
``sequentially'', the remaining converses have been shown to fail, which has
been summarised in \cite[Proposition 3.1]{KD1}$.$ We will return to some of
these nonreversible implications in Section 6.

\section*{{\bf 4}. A Review of Forms of Weak Supercyclicity}

\vskip6pt
There is no reason to investigate any form of weak cyclicity since, as we saw
above, weak cyclicity coincides with (strong) cyclicity$.$ Our focus in this
section is narrowed exclusively to
\vskip4pt\noi
\centerline{forms of\, {\it weak supercyclicity}\/ \,only,}
\vskip6pt\noi
as described in Section 3$.$ Historically, weak supercyclicity followed the
first works on weak hypercyclicity$.$ So, before discussing forms of weak
supercyclicity properly, we take as a starting point a few pioneering works
on weak hypercyclicity that ex\-plicitly motivated weak
supercyclicity$.\kern-1pt$ Since the purpose of this section is to present a
survey on forms of weak supercyclicity, aiming to discuss the literature in a
rough chronological order, all references cited in this section will be
followed by their respective year of publication.

\vskip6pt\noi
{\bf \small 4.1.}
{\sc pre-weak supercyclicity$:$ 2002--2004}$.$
Feldman \cite{Fel1} (2002) asked whether every weakly hypercyclic operator
is hypercyclic (i.e., strong hypercyclic --- in the norm topology)$.$ Chan and
Sanders \cite{CS} (2004) answered that question by exhib\-iting a weakly
hypercyclic operator that is not hypercyclic, and they have asked if the
spectrum of a weakly hypercyclic operator meets the unit circle$.$ As it
happened in the norm topology case shown by Kitai \cite{Kit}
(1982)$\kern.5pt$), Dilworth and Troitsky \cite{DT} (2003) have shown that
every component of the spectrum of a weakly hypercyclic operator intersects
the unit circle, thus answering the question posed in \cite{CS} (2004)$.$
Although the above mentioned works did not deal with weak supercyclicity, but
with hypercyclicity, weak or strong, they motivated and influenced the
subsequent works on the several forms of weak supercyclicity.

\vskip6pt\noi
{\bf \small 4.2.}
{\sc the beginning of weak supercyclicity$:$ 2004--2005}$.$
It seems that the first results specifically on weak supercyclicity are due
to Sanders in \cite{San1} (2004) and in \cite{San2} (2005)$.$ In \cite{San1}
(2004) she constructed a weakly supercyclic counterpart of the above mentioned
result from \cite{CS} (2004) by exhibiting a weakly
\hbox{supercyclic operator} that is not weakly hypercyclic$.$ Also, it was
exhibited in \cite{San1} (2004) a weakly super\-cyclic operator that is not
supercyclic, and it was also shown that the set of weakly supercyclic vectors
is dense in the norm topology$.$ Moreover, it was shown in \cite{San1} (2004)
that hyponormal operators on a Hilbert space are not weakly hypercyclic,
neither supercyclic as shown by Bourdon \cite{Bou} (1997)$\kern.5pt$, and the
question of whether there exists a weakly supercyclic hyponormal operator
closed the paper \cite{San1} (2004) (which was affirmatively answered later,
as we will see below)$.$ In contrast with a result by Ansari and Bourdon
\cite{AB} (1997), which says that there is no supercyclic isometry on Banach
spaces, it was exhibited in \cite{San2} (2005) a weakly supercyclic isometry
on the Banach space $c_0(\ZZ)$ of bilateral sequences converging to zero on
both sides$.$ In \cite{Pra} (2005) Pr\v ajitur\v a gave a collection of
spectral properties for weakly supercyclic operators on a Hilbert space $\H.$
(For supercyclic operators, a spectral description had been given by Herrero
in \cite{Her} (1991).) It was also shown in \hbox{\cite{Pra} (2005)} that the
collection of all operators that are not weakly supercyclic is dense in $\BH$.

\vskip6pt\noi
{\bf \small 4.3.}
{\sc further forms of weak supercyclicity$:$ 2005--2006}$.$
Although the term weak l-sequential supercyclicity had not been coined at that
time, B\`es, Chan, and Sanders considered the notion in \cite{BCS1} (2005) to
show that a bilateral weighted shift, on $\ell^p$ for ${1\le p<\infty}$ or
on $c_0(\ZZ)$, is weakly l-sequentially supercyclic if and only if it is
(strong) supercyclic$.$ Forms of weak* sequential supercyclicity have also
been investigated by B\`es, Chan, and Sanders in another paper \cite{BCS2}
(2006)$.\kern-1pt$ The notion of weak l-sequential supercyclicity was also
implicitly used in at least one part (proof of Lemma 3.3) of \cite{BM1}
(2006)$.$ Bayart and Matheron \cite{BM1} (2006) proved that if a hypo\-normal
operator is weakly supercyclic, then it is a multiple of a unitary$\kern-1pt$
They also proved the existence of a weakly supercyclic unitary operator, in
fact, the existence of a weakly \hbox{l-sequentially} supercyclic unitary
operator, thus answering the question posed in \cite{San1} (2004) --- there
are weakly supercyclic hyponormal operators --- and also extending to Hilbert
spaces the existence of weakly supercyclic isometries on $c_0(\ZZ)$ proved in
\cite{San2} (2005)$.\kern-1pt$ The question of whether there exists a weakly
hyper\-cyclic operator that is not weakly sequentially supercyclic closed the
paper \cite{BM1} (2006) (which was affirmatively answered later, as we will
see below).

\vskip6pt\noi
{\bf \small 4.4.}
{\sc clarifying three forms of weak supercyclicity$:$ 2007}$.$
The notions of weak l-sequential supercyclicity, weak sequential
supercyclicity, and weak supercyclicity, as presented in Section 3, have been
considered by Shkarin \cite{Shk1} (2007), where the notion of weak
l-sequential supercyclicity, stronger than the usual weak sequential
supercyclicity but weaker than (strong) supercyclicity, was discussed in
detail, thus amounting to the following inclusions:
\vskip5pt\noi
\centerline{\esc
weak l-sequential supercyclicity
$\,\scriptstyle\sse\,$
weak sequential supercyclicity}
\centerline{\esc
$\,\scriptstyle\subset\,$
weak supercyclicity.}
\vskip6pt\noi
Special attention was given to whether these were proper inclusions$.$ It was
shown in \cite{Shk1} (2007) the existence of a weakly supercyclic unitary
operator that is not weakly sequentially supercyclic, and so not weakly
l-sequentially supercyclic$.$ It was also asked in \cite{Shk1} (2007) whether
there is a weakly sequentially supercyclic operator that is not weakly
l-sequentially supercyclic$.$ (Is there?) $\kern-1pt$The original terminology
for weak l-sequential supercyclicity is different from ours: we use the letter
``l'' for ``limit'', borrowed from \cite{KD1} (2018), while the numeral ``1''
was used in \cite{Shk1} (2007) --- and there are reasons for both
terminologies.

\vskip6pt\noi
{\bf \small 4.5.}
{\sc weak supercyclicity for operators close to normal$:$ 2006--2020}$.$
Consider the well-known chain of classes of Hilbert-space operators$:$
\vskip4pt\noi
\centerline{\esc
unitary
$\,\scriptstyle\subset\,$
normal
$\,\scriptstyle\subset\,$
hyponormal
$\,\scriptstyle\subset\,$
quasihyponormal}
\centerline{\esc
$\,\scriptstyle\subset\,$
semi-quasihyponormal
$\,\scriptstyle\subset\,$
paranormal
$\,\scriptstyle\subset\,$
normaloid.}
\vskip6pt\noi
For the above classes, set $|T|=(T^*T)^{1/2}\!$ and recall that an operator
$T\kern-.5pt$ on \hbox{a Hilbert} space $\H$ is, by definition, unitary if
${|T^*|^2=|T|^2=I}$, normal if ${|T^*|^2=|T|^2}\kern-.5pt$, hypo\-normal if
${|T^*|^2\kern-.5pt\le|T|^2}\kern-.5pt$, quasihyponormal if
${|T|^4\kern-.5pt\le|T^2|^2}\kern-.5pt$, semi-quasihyponormal if
${|T|^2\kern-.5pt\le|T^2|}$ (also known as class $\A$), paranormal if
${\|Tx\|^2\kern-.5pt\le\|T^2x\|\,\|x\|}$ for every $x$ in $\H$, and normaloid
if ${\|T^n\|=\|T\|^n}$ for every ${n\ge0}$, equivalently, if ${r(T)=\|T\|}$.

\vskip6pt\noi
There is no supercyclic unitary operator$.$ Indeed, as we have seen in
Subsection 4.2, Ansari and Bourdon \cite[Proof of Theorem 2.1]{AB} (1997)
have shown that there is no supercyclic isometry on a Banach space --- as a
matter of fact, Hilden and Wallen \cite[p.564]{HW} (1974) had already shown
that the unilateral shift is not supercyclic$.$ Actually, Bourdon
\cite[Theorem 3.1]{Bou} (1997) had shown that there is no supercyclic
hypo\-normal operator on a Hilbert space, as we have also seen in Subsection
4.2, extending Kitai's result \cite[Corollary 4.5]{Kit} (1982) where she had
shown that there is no hypercyclic hyponormal operator; and also extending
Hilden and Wallen's result \cite[p.564]{HW} (1974) where they had shown that
there is no supercyclic normal operator$.$ In spite of this, there are
weakly supercyclic unitary (thus normal, and hence hyponormal) operators on a
Hilbert space$.$ In fact, we saw in Subsection 4.3 that Bayart and Matheron
\cite[Example 3.6]{BM1} (2006) have shown that there exist weakly
l-sequentially supercyclic (thus weakly supercyclic) unitary operators, and
also that if a hyponormal operator is weakly supercyclic, then it is a
multiple of a unitary operator \cite[Theorem 3.4]{BM1} (2006)$.$ Shkarin
\cite[Corollary 1.3]{Shk1} (2007) has shown the existence of a weakly
supercyclic unitary operator that is not weakly sequentially (thus not weakly
l-sequentially) supercyclic$.$ Results about reduction to multiples of unitary
operators have been extended beyond the hyponormal op\-erators$.$ For
instance, Jung, Kim, and Ko \cite[Theorem 4.2]{JKK} (2014) proved that if an
injective semi-quasihyponormal (i.e., class $\A$) operator is weakly
supercyclic, then it is a scalar multiple of a unitary$.$ This was followed by
Duggal, Kubrusly, and Kim \cite[Theorem 2.7]{DKK} (2015), where the
injectivity assumption was dismissed$.$ Duggal \cite[Corollary 3.2]{Dug}
(2016) has shown that if $T$ is such that
${\|T^*x\|^2\kern-1pt\le\kern-.5pt\|T^2x\|\kern.5pt\|x\|}$ for \hbox{every}
${x\kern-1pt\in\kern-1pt\H}$ (such operators are called $*$-paranormal) and is
weakly supercyclic, then it is a scalar multiple of a unitary (among other
results)$.$ Still along this line, Shen and Ji \cite{SJ} (2020) have shown
that weakly supercyclic quasi-2-expansive operator (i.e., an operator $T$ for
which
${T^{*3}T^3\!-\kern-1pt2\kern1ptT^{*2}T^2\!+\kern-1ptT^*T
\hskip-1pt\le\hskip-1pt O}$)
is unitary$.\kern-1pt$ \hbox{All the above} classes include the unitary
operators, as discussed in \cite{BM1}, \cite{Shk1}, \cite{JKK}, \cite{DKK},
\cite{Dug}, \cite{SJ}$.$ Weak forms of supercyclicity for unitary operators
will be considered in Section 6.

\vskip6pt\noi
{\bf \small 4.6.}
{\sc weak supercyclicity for compact operators$:$ 2007--2018}$.$
It had already been shown by Le\'on-Saavedra and Piqueiras-Lerena \cite{LP}
(2003) and Gallardo-Guti\'errez and Montes-Rodr\'\i ques \cite{GM} (2004) that
the Volterra operator on $L^p[0,1]$, ${p\ge1}$ is not supercyclic$.\kern-1pt$
A nonelementary proof that the Volterra operator is not even weakly
supercyclic was given by Montes-Rodr\'\i gues and Shkarin \cite{MS} (2007),
and so, in particular, it is not weakly l-sequentially supercyclic$.$ Shkarin
\cite{Shk2} (2010) has shown that any operator on $L^p[0,1]$ commuting with
the Volterra operator is not weakly supercyclic$.$ Also, in a paper linked to
Volterra operators, Hedayatian and Faghih-Ahmadi \cite{HF} have observed that
there is no weakly supercyclic \hbox{operator} in the commutant of a cyclic
convolution operator$.$ It has been verified by Kubrusly and Duggal \cite{KD2}
(2018) that for a compact operator on a normed space, weak l-sequential
supercyclicity coincides with supercyclic$.\kern-1pt$ Thus, since Volterra
operators are compact and not supercyclic, as we saw above, the non-weak
l-sequential supercyclicity comes at once$.$ It was also shown in
\cite{KD2} (2018) that a weakly l-sequentially supercyclic compact operator 
on an infinite-dimensional Banach space is quasinilpotent (recall that
Volterra operators are also quasinilpotent)$.$

\vskip6pt\noi
{\bf Remark 3.}$\kern-.5pt$ {\sc still on compact operators}\/$.$
\vskip3pt\noi
(i)
As weak l-sequential supercyclicity coincides with supercyclicity for compact
operators, weak supercyclicity coincides with supercyclicity for bilateral
weighted shifts on $\ell^p$ for ${1\le p<2}$ \cite[Theorem 6.3]{MS} (2007) --- 
for weak l-sequential supercyclicity and for ${1\le p<\infty}$, see
\cite[Theorem 1]{BCS1} (2005), as mentioned in 4.3 above.
\vskip4pt\noi
(ii)
A hypercyclic operator is not compact, and there is no supercyclic operator
on a complex normed space with a finite dimension greater than 1
\cite[\hbox{Section 4}]{Hez} (1992)$.$ There are, however, compact supercyclic
operators on a separable infinite-dimensional complex Banach space
\cite[Theorem 1 and Section 4]{Hez}$.$ Consequently, there are weakly
supercyclic compact operators for any form of weak supercyclicity.
\vskip4pt\noi
(iii)
The Volterra operator is an example of a compact quasinilpotent that is not
super\-cyclic, and so it is not weakly l-sequentially supercyclic, as we saw
above$.$ But, as we saw in item (b), there exist quasinilpotent supercyclic
operators$.$ It \hbox{was shown} in \cite[Corollary 5.3]{Sal} (1999) that if
the adjoint of a bilateral weighted shift on $\ell^2\!$, or of a unilateral
weighted shift on $\ell_+^2$, has a weighting sequence possessing a
subsequence that goes to zero, then there is an infinite-dimensional
subspace whose all nonzero vectors are supercyclic for it$.$ In particular,
this holds for the adjoint of weighted shifts with weighting sequences
converging to zero, and so this happens for the adjoint of compact
quasinilpotent weighted shifts.

\vskip6pt\noi
{\bf \small 4.7.}
{\sc weak supercyclicity for composition operators$:$ 2010--2021}$.$
There are a large number of works dealing with cyclicity and hypercyclicity
(weak or not), as well as with supercyclicity, for composition operators$.$
The number, however, is not so large when dealing with weak supercyclicity$.$
Kamali, \hbox{Hedayatian}, and Khani Robati \cite{KHK} (2010) investigated
conditions for weak supercyclicity of weighted composition operators, in
particular on the Hilbert space of holomorphic functions on the unit
disk$.\kern-1pt$ Weak supercyclicity, and weak hypercyclicity, of composition
operators on some Hilbert spaces of analytic functions (e.g., on certain
weighted Hardy spaces) were investigated by Kamali \cite{Kam} (2013)$.$
B\`es (2013) has given necessary and sufficient conditions for a composition
operator to be weakly supercyclic$.$ Conditions for weak supercyclicity of
some classes of weighted composition operators were given by Moradi, Khani
Robati, and Hedayatian \cite{MKH} (2017), where specific spaces on which weak
supercyclicity of composition operators never holds were also investigated$.$
Beltr\'an-Meneu, Jord\'a, and Murillo-Arcila \cite{BJM} (2020) also considered
conditions ensuring that certain weighted composition operators are not weakly
supercyclic$.$ See also Mengestie \cite{Men} (2021) for non-weak
supercyclicity of weighted composition operators on certain Fock spaces.

\vskip5pt\noi
{\bf \small 4.8.}
{\sc weak supercyclicity and weak stability$:$ 2016--2024}$.$
This is a central topic in the present exposition$.$ A review of current
literature dealing with weak supercyclicity and its relation to weak
stability will be given in Section 6$.$ Before discussing if (or when) weak
supercyclicity implies weak stability, we need to consider the motivating
fact that is behind such a problem, viz., supercyclicity implies strong
stability$.\kern-1pt$ This fundamental result will be treated in detail in
the \hbox{next section}.

\vskip5pt\noi
{\bf Remark 4.} $\kern-.5pt$ {\sc further results related to weak
supercyclicity}\/$.$
\vskip3.5pt\noi
(i)
The term {\it weak supercyclicity}\/ has also been used to mean a notion
different from those considered above$.$ For instance, a subset $S$ of a
Hilbert space $\H$ is said to be {\it $n$-weakly dense}\/ in $\H$ if for every
$n$-dimensional subspace $\M$ of $\H$, the image of $S$ under the orthogonal
projection $E$ onto $\M$ is dense in $\M$ (i.e., if ${E(S)^-\!=\M}$ for the
orthogonal projection ${E\!:\H\to\H}$ such that ${\rm{range}\kern1pt(E)=\M}$
if $\dim\M=n).$ An operator $T$ in $\BH$ is {\it $n$-weakly supercyclic}\/ if
it has a projective orbit $\Oe([y])$ for some $y$ in $\H$ that is $n$-weakly
dense in $\H$ (see, e.g., \hbox{\cite{Fel2,Fel3} (2011, 2012)}$\kern.5pt$).

\vskip3.5pt\noi
(ii)
If an operator is invertible and supercyclic, then so is its inverse
\cite[Section 4]{AB} (1997), \cite[Corollary 2.4]{Sal} (1999)$.$ There are,
however, invertible weakly supercyclic operators on $\ell^p$ for
${p\in[2,\infty)}$ whose inverses are not weakly supercyclic
\cite[Corollary 2.5]{San1} (2004)$.$ Is the inverse (i.e., the adjoint) of a
weakly supercyclic unitary operator weakly supercyclic$?$ Recall that the
adjoint of a supercyclic coisometry may not be supercyclic$.$ Example$:$ a
backward unilateral shift $S^*$ is a supercyclic coisometry
\cite[Theorem 3]{HW} (1974), while its adjoint, the unilateral shift $S$,
being an isometry, is not supercyclic \cite[p.564]{HW} (1974) --- also see
\cite[Proof of Theorem 2.1]{AB} (1997) or \cite[Lemma 4.1(b)]{KD1} (2018)$.$
The same happens with weak supercyclicity$:$ the adjoint of a weakly
supercyclic coisometry may not be weakly supercyclic$.$ Example$:$ $S$ is not
weakly supercyclic --- since a weakly supercyclic hyponormal operator is a
multiple of a unitary operator \cite[Theorem 3.4]{BM1} (2006) --- but
$S^*\kern-1pt$ is weakly supercyclic, since it is supercyclic$.$ However, if
an Hilbert-space isometry is invertible (i.e., if it is unitary) and weakly
l-sequentially supercyclic, then it has a weakly l-sequentially supercyclic
adjoint \cite[Theorem 3.3]{KD2} (2018).
\vskip3.5pt\noi
(iii)
The notion of hypercyclicity, in particular, has been investigated in
Fr\'echet spaces, or F-spaces, or locally convex spaces (see, e.g., \cite{Ans}
(1997), \cite{BP} (1998), \cite{Per} (2001), \cite{BFPW} (2005), \cite{BM2}
(2009), \cite{GP} (2011), among others)$.$ Our discussion in the present paper
is restricted to normed (or Banach) spaces.
\goodbreak

\section*{{\bf 5}. Strong Supercyclicity and Strong Stability}

An important motivation for an investigation on weak supercyclicity and its
relationship with weak stability is the fact that (strong) supercyclicity
and strong stability behave nicely$:$ {\it supercyclicity implies strong
stability}\/$.$ The verification of this fact is all but trivial, and the
proof itself naturally impels the quest for a weak counterpart$.$ For this
reason, we begin by looking at the strong case first.

\vskip6pt
It was proved in \cite[Theorem 2.1]{AB} that {\it if a power bounded
operator\/ $T\!$ \hbox{on a Banach} space\/ $\X$ is such that\/
$\|T^nx\|\not\to0$ for every nonzero\/ ${x\in\X}$, then\/ $T\!$ has no
supercyclic vector}\/$.$ In other words, a power bounded operator of class
$C_{1{\textstyle\cdot}}$ is not supercyclic$.$ Based on this result, more was
proved in \cite[Theorem 2.2]{AB}, namely, {\it if a power bounded operator\/
${T\kern-1pt\in\kern-1pt\BX}$ is such that there exists a vector\/ ${y\in\X}$
for which there is a scalar sequence\/ $\{\alpha_k\}_{k\ge0}$ and a sequence\/
$\{T^{n_k}\}_{k\ge0}$ with entries from\/ $\{T^n\}_{n\ge0}$ such that
${\alpha_kT^{n_k}y\to x}$ for each\/ ${x\in\X}$, then\/ ${T^nz\to0}$ for
every}\/ ${z\in\X}.$ In other words, a power bounded supercyclic operator on a
normed space is \hbox{strongly stable}.

\vskip6pt
These significant results by Ansari and Bourdon \cite{AB} (1997) are essential
in linking supercyclicity to stability, establishing a direct implication from
(strong) supercyclicity to strong stability$.$ Such an implication will be
crucial in any attempt to extend it, expecting that some form of weak
supercyclicity will imply weak stability$.$

\vskip6pt
The above mentioned results are stated below, followed by detailed proofs$.$
This is done here not only for sake of completeness but mainly because a
close analysis of the structure of their proofs will be fundamental for the
next section$.$ In fact, the arguments in the next proofs are required to
yield a possible counterpart of Theorem 5.2 below for weak l-sequential
supercyclicity and weak stability.

\vskip6pt\noi
{\bf Theorem 5.1.}$\,$\cite{AB}
{\it A power bounded operator of class\/ $C_{1{\textstyle\cdot}}\kern-1pt$
acting on a normed space is not super\-cyclic}\/.

\vskip6pt\noi
{\it Proof}\/$.$
The well-known notion of Banach limit is a central argument in the proof$.$
For a summary of all Banach limit properties required here, see, e.g.,
\cite[\hbox{Section 6}]{KD3}$.$ Let the functional
${\vphi\!:\ell_+^\infty\!\to\CC}$, assigning a complex number
$\vphi(\{\xi_n\})$ to each bounded complex-valued sequence
$\{\xi_n\}={\{\xi_n\}_{n\ge0}\in\ell_+^\infty}$, be a Banach
limit$.\kern-1pt$ This \hbox{means that} $\vphi$ is linear, bounded, positive
(i.e., ${0\le\vphi(\{\xi_n\})}$ whenever ${0\le\xi_n}$ for every integer
${n\ge0}$), backward-shift-invariant (i.e.,
${\vphi(\{\xi_n\})=\vphi(\{\xi_{n+1}\})}\kern.5pt$), with ${\|\vphi\|=1}$,
and assigns to each convergent sequence its own limit$.$ Let
$({\X,\|\cdot\|})$ be any \hbox{normed space}$.$ Suppose ${T\in\BX}$ is power
bounded$.$ For each ${x\in\X}$ consider the bounded sequence of nonnegative
numbers $\{\|T^nx\|\}_{n\ge0}$, and set
$$
\|x\|_\vphi\!=\vphi(\{\|T^nx\|\})
$$
for any Banach limit $\vphi.$ Since a Banach limit $\vphi$ is
order-preserving for real-valued bounded sequences (i.e.,
$\vphi(\{\xi_n\})\le\vphi(\{\upsilon_n\})$ whenever $\xi_n$ and $\upsilon_n$
are real numbers such that ${\xi_n\le\upsilon_n}$), with
$\vphi(\{1,1,1,\dots\})=1$, and since $T$ is power bounded, this defines an
$\RR$-valued functional ${\|\cdot\|_\vphi\!:\X\!\to\RR}$ with the property
that, for \hbox{every ${x\in\X}$},
$$
\|x\|_\vphi\!\le({\sup}_n\|T^n\|)\,\|x\|.
$$
But this defines a seminorm on $\X.$ Indeed, since $\vphi$ is a positive
functional, we get ${0\le\|x\|_\vphi}$ for every ${x\in\X}.$ Since, in
addition, $\vphi$ is linear (as well as $T$) and order-preserving, we get
$|\alpha|\|x\|_\vphi=\|\alpha x\|_\vphi$ and
$\|{x+y}\|_\vphi\le\|x\|_\vphi\!+\|y\|_\vphi$ for every ${x,y\in\X}$ and
every ${\alpha\in\CC}.$ However, as the power bounded operator $T$ is of
class $C_{1{\textstyle\cdot}}$, it follows that ${\|\cdot\|_\vphi}$ is a
norm, as we will see next.

\vskip6pt\noi
{\it Claim 0}\/$.$
If $T$ is power bounded operator on a normed space, then
$$
T^nx\not\to0
\quad\;\hbox{implies}\;\quad
0<{\liminf}_n\|T^nx\|
\quad\;\hbox{which implies}\;\quad
0<{\inf}_n\|T^nx\|.
$$
\vskip-4pt

\vskip6pt\noi
{\it Proof of Claim 0}\/$.$
(i)
Suppose $\{\beta_n\}_{n\ge0}$ is a real-valued sequence such that
$$
0\le\beta_n\le\beta\,\beta_m
\quad\;\hbox{whenever}\;\quad
m\le n
$$
for some positive constant $\beta.$ If ${\beta_n\not\to0}$, then
${0<\beta_n\le\beta\,\beta_0}$ for every $n.$ So $\{\beta_n\}_{n\ge0}$ is a
bounded sequence of positive numbers for which there is a subsequence
$\{\beta_{n_i}\}_{i\ge0}$ such that ${0<\lim_i\beta_{n_i}}$ and
${\beta_{n_i}\le\beta\,\beta_n}$ if ${n<n_i}.$ Hence
${0<\lim_i\beta_{n_i}}={\liminf_i\beta_{n_i}}\le{\beta\,\liminf_n\beta_n}$,
which ensures ${0<\inf_n\beta_n}$ (since ${0<\beta_n}$ for every $n$)$.$ Thus
$$
\beta_n\not\to0
\qquad\limply\qquad
0<{\liminf}_n\beta_n
\qquad\limply\qquad
0<{\inf}_n\beta_n.
$$
\vskip0pt\noi
(ii)
Since $0\le\|T^nx\|=\|{T^{n-m}\,T^m}x\|\le(\sup_n{\|T^n\|)\,\|T^mx\|}$
for every $x$ whenever ${m\kern-1pt\le\kern-1pt n}$, the proof of Claim 0
is complete by setting ${\beta_n\!=\kern-1pt\|T^nx\|}$ and
${\beta\kern-1pt=\kern-1pt\sup_n\kern-1pt\|T^n\|}.\!\!\!\qed$

\vskip6pt\noi
So, according to Claim 0, if $T$ is of class $C_{1{\textstyle\cdot}}\kern-1pt$
(i.e., if ${T^nx\not\to0}$ for every ${0\ne x\in\X}$), then
${0<\liminf_n\|T^nx\|}$ for every ${0\ne x\in\X}.$ Moreover, as
$\liminf_n\xi_n\le\vphi(\{\xi_n\})\le\limsup_n\xi_n$ for every real-valued
bounded sequence $\{\xi_n\}$, we get $0<\vphi(\{\|T^nx\|\})=\|x\|_\vphi$
if ${x\ne0}.$ Consequently,
$$
\hbox{$T\kern-1pt$ is of class $C_{1{\textstyle\cdot}}
\quad\limply\quad\|\cdot\|_\vphi$ is a norm on $\kern-1pt\X$.}
$$
Thus assuming that $T\kern-1pt$ is of class $C_{1{\textstyle\cdot}}\kern-1pt$
is required to ensure that $\|\cdot\|_\vphi$ is a norm$.$ Then consider the
normed space ${\X_\vphi=(\X,\|\cdot\|_\vphi)}$ wherein the operator $T$ acts
as an isometry$.$ In fact, since $\vphi$ is backward-shift-invariant, we get
for every ${x\in\X}$
$$
\|Tx\|_\vphi=\vphi(\{\|(T^{n+1}x)\|\})
=\vphi(\{\|(T^nx)\|\})=\|x\|_\vphi.
$$
Now consider the completion ${\what\X\!=(\what\X,\|\cdot\|_{\what\X})}$ of
the normed space ${\X_\vphi\!=(\X,\|\cdot\|_\vphi)}$, so tat $\X_\vphi$ is
isometrically and densely embedded in $\what\X$, which means $\X_\vphi$ is
isometrically isomorphic to a dense linear manifold $\wtil\X$ of the Banach
space $\what\X.$ In other words, there is an isometric isomorphism
${J\!:\X_\vphi\kern-1pt\to J(\X_\vphi)=\wtil\X}$ where $\wtil\X$ is dense in
$\what\X.$ Consider the extension over completion
${\what T\kern-1pt\in\B[\what\X]}$ of ${T\kern-1pt\in\kern-1pt\B[\X_\vphi]}$,
and let $\wtil T=$ ${J\kern1ptTJ^{-1}}$ be the restriction of $\what T$ to
$\wtil\X.$ That is,
$\wtil T=\what T|_{\wtil\X}={J\kern1ptTJ^{-1}\!\in\kern-1pt\B[\wtil\X]}.$ Thus
$\wtil\X$ is $\what T$-invariant and
$T
=J^{-1}\wtil TJ
=J^{-1}\what T|_{\wtil\X}J
={J^{-1}\what TJ\in\kern-1pt\B[\X_\vphi]}$,
as in the following commutative diagram, where the symbol
$\kern-1pt\dsse\kern-1pt$ means densely included:
\vskip5pt\noi
$$
\matrix{
\kern7pt\X_\vphi\phantom{\int_|}\kern-5pt
& \kern-10pt\buildrel J\over\emap
& \kern-20pt\wtil\X\;=\;J(\X_\vphi)
& \kern-35pt\dsse
& \kern-25pt\what\X                                                   \cr
\kern-5pt\scriptstyle{T}\kern1pt\nmap
&
& \kern-7pt{\nmap\scriptstyle\wtil T=\what T|_{\wtil\X}\,=J\,TJ^{-1}}
&
& \kern-20pt\kern1pt\nmap\scriptstyle{\what T}                        \cr
\kern7pt\X_\vphi
& \kern-10pt\buildrel J\over\emap
& \kern-20pt \wtil\X\;=\;J(\X_\vphi)
& \kern-35pt\dsse
& \kern-25pt\;\what\X.                                                \cr}
$$
\vskip5pt\noi
Next suppose $T$ is supercyclic, and observe that supercyclicity easily
survives the change of norms from ${\|\cdot\|}$ to ${\|\cdot\|_\vphi}$;
that is,
\vskip3pt\noi
\begin{description}
\item{$\kern-6pt$(a)$\kern3pt$}
\it if ${y\in\X}$ is a supercyclic vector for $T$ when it acts on
${\X=(\X,\|\cdot\|})$, then this same ${y\in\X}$ is a supercyclic vector for
$T$ when it acts on ${\X_\vphi=(\X,\|\cdot\|}_\vphi)$.
\end{description}
\vskip3pt\noi
Indeed, since ${0<{\sup}_n\|T^n\|<\infty}$ and
$\|z\|_\vphi\!\le({\sup}_n\|T^n\|)\,\|z\|$ for every ${z\in\X}$, it follows
at once that
${\|{\alpha_iT^{n_i}y-x}\|_\vphi\le(\sup_n\|T^n\|)\|{\alpha_iT^{n_i}y-x}\|}$,
and so $y$ is a supercyclic vector in $\X_\vphi$ whenever it is a supercyclic
vector in $\X$.
\vskip6pt\noi
This ensures that supercyclicity can be extended from $T$ on $\X$ to $\what T$
on $\what\X.$ In fact,
\vskip3pt\noi
\begin{description}
\item{$\kern-6pt$(b)$\kern3pt$}
\it if\/ ${y\in\X}$ is a supercyclic vector for\/ $T\!$ when it acts on\/
$\X_\vphi$, then its image\/ $\wtil y=J{y\in\wtil\X}$ under the isometric
isomorphism\/ $J$ of\/ $\X_\vphi$ onto the dense linear manifold\/ $\wtil\X$
of the Banach space\/ $\what\X$ is supercyclic for the
\hbox{extension\/ $\what T\!$ on\/ $\what\X\!$ of\/ $T\!$.}
\end{description}
\vskip3pt\noi
The above result was applied in the proof of Theorem 2.1 from \cite{AB} (on
p.198, lines 11,12) but a proof was omitted$.$ We give a proof of it in
Claim I below.

\vskip6pt\noi
{\it Claim I}\/$.$
If there is a supercyclic vector ${y\in\X_\vphi}$ for
${T\kern-1pt\in\kern-1pt\B[\X_\vphi]}$, then there is a supercyclic vector
${\wtil y\in\what\X}$ for ${\what T\kern-1pt\in\B[\what \X]}$.

\vskip6pt\noi
{\it Proof of Claim I}\/$.$
Since
${\wtil T
\kern-1pt=\kern-1pt
\what T|_{\wtil\X}
\kern-1pt=\kern-1pt
{J\kern1pt T\kern-.5ptJ^{-1}}}\kern-1pt$,
where $\wtil\X$ is $\what T$-invariant, we get 
$\wtil T^n
\kern-1.5pt=\kern-1pt
{(\what T|_{\wtil\X})^n
\kern-1pt=\kern-1pt
\what T^n|_{\wtil\X}
\kern-1pt=\kern-1.5pt
J\kern1ptT^n\kern-1ptJ^{-1}}\kern-1pt$
for every integer ${n\kern-1pT\ge\kern-1pt0}.$ Now split the proof into
\hbox{two parts}.

\vskip6pt\noi
(1)
If there is a supercyclic vector ${y\in\X_\vphi}$ for
${T\kern-1pt\in\kern-1pt\B[\X_\vphi]}$, then there is a supercyclic vector
${\wtil y\in\wtil\X}$ for ${\wtil T\kern-1pt\in\kern-1pt\B[\wtil\X]}$.

\vskip6pt\noi
(2)
If there is a supercyclic vector ${\wtil y\in\wtil\X}$ for
${\wtil T\kern-1pt\in\B[\wtil\X]}$, then $\wtil y$ is a supercyclic vector for
${\what T\kern-1pt\in\B[\what \X]}$.

\vskip6pt\noi
{\it Proof of\/ \rm(1)}\/$.$
If there is a supercyclic vector ${y\in\kern-.5pt\X_\vphi}\kern-.5pt$ for
${T\kern-1pt\in\kern-.5pt\B[\X_\vphi]}$, then for \hbox{every}
${0\ne x\in\kern-.5pt\X_\vphi}\kern-.5pt$ there exists a sequence of nonzero
numbers $\{\alpha_j\}_{j\ge0}\kern-.5pt$ that depends on $x$, and so we write
$\{\alpha_j(x)\}_{j\ge0}$, and a sequence $\{T^{n_j}\}_{j\ge0}\kern-.5pt$
with entries \hbox{from $\{T^n\}_{n\ge0}$}, which also depends on $x$, such
that, in ${\X_\vphi\kern-.5pt=(\X,\|\cdot\|_\vphi)}$,
$$
\alpha_j(x)\kern1ptT^{n_j}y\conv x.
$$
Set ${\wtil y=Jy}$ in $\wtil\X.$ As every $\wtil x$ in $\wtil\X$ is of the
form $\wtil x=Jx$ for some $x$ in $\X_\vphi$, thus for every $\wtil x$ in
$\wtil\X$ set $\wtil\alpha_j(\wtil x)=\wtil\alpha_j(Jx)=\alpha_j(x)$ for some
${x\in\X_\vphi}.$ Take an arbitrary ${\wtil x\in\wtil\X}.$ Then there exists a
sequence $\{\wtil\alpha_j(\wtil x)\}_{j\ge0}$ such that
(with ${\wtil T^{n_j}\kern-1pt=J\,T^{n_j}J^{-1}}$)
\begin{eqnarray*}
\|\wtil\alpha_j(\wtil x)\kern1pt\wtil T^{n_j}\wtil y-\wtil x\|_{\what\X}
&\kern-6pt=\kern-6pt&
\|\alpha_j(x)\kern1ptJ\kern1ptT^{n_j}y-Jx\|_{\what\X}                     \\
&\kern-6pt=\kern-6pt&
\|J(\alpha_j(x)\kern1ptT^{n_j}y-x)\|_{\what\X}
=\|\alpha_j(x)T^{n_j}y-x\|_{\X_\vphi}
\end{eqnarray*}
for some ${x\in\X_\vphi}.$ Hence, since ${y\in\X_\vphi}$ is a supercyclic
vector for ${T\in\B[\X_\vphi]}$ such that
${\|\alpha_j(x)T^{n_j}y-x\|_{\X_\vphi}\!\to0}$, we get
${\|\wtil\alpha_j(\wtil x)\kern1pt\wtil T^{n_j}\wtil y-\wtil x\|_{\what\X}\!
\to0}$, and therefore $\wtil y={Jy\in\wtil\X}$ is a supercyclic vector
for ${\wtil T\in\B[\wtil\X]}.$ This proves (1): for every
${\wtil x\in\wtil\X}$
$$
\wtil\alpha_j\kern1pt(\wtil x)\kern1pt\wtil T^{n_j}\wtil y\conv \wtil x.
$$
\vskip-2pt

\vskip6pt\noi
{\it Proof of\/ \rm(2)}\/$.$
Take an arbitrary ${\what x\in\what\X}.$ Since ${\wtil\X^-\!=\what\X}$, there
is a sequence $\{\wtil x_k\}_{k\ge0}$, with
${\wtil x_k=\wtil x_k(\what x)\in\wtil\X}$ for each $k$, such that
${\wtil x_k\to\what x}.$ Then, by the above convergence,
$$
\wtil\alpha_j\kern1pt(\wtil x_k)\kern1pt\wtil T^{n_j}\wtil y
\jconv\wtil x_k\kconv\what x.                                       \eqno(*)
$$
This ensures the existence of a sequence of nonzero numbers
$\{\what\alpha_i(\what x)\}_{i\ge1}$ such that
$$
\what\alpha_i(\what x)\kern1pt\what T^{n_i}\wtil y\conv\what x     \eqno(**)
$$
for a sequence\/ $\{\what T^{n_i}\}_{i\ge1}$ with entries from
$\{\what T^n\}_{n\ge0}$, where each $\what T^n$ is the extension over
completion of each $T^n$ --- equivalently, the extension by continuity of
each $\wtil T^n\kern-1pt$.
\vskip6pt\noi
Indeed, consider both convergences in ($*$)$.$ Take an arbitrary ${\veps>0}.$
Thus there exists a positive integer $k_\veps$ such that
${\|\wtil x_k-\what x\|_{\what\X}\le\frac{\veps}{2}}$ whenever
${k\ge k_\veps}.$ Moreover, for each $k$ there exists a positive integer
$j_{\kern1pt\veps,k}$ such that
${\|\wtil\alpha_j(\wtil x_k)\kern1pt\wtil T^{n_j}\wtil y-\wtil x_k\|_{\what\X}
\le\frac{\veps}{2}}$
whenever ${j\ge j_{\veps,k}}.$ Therefore, by taking ${k=k_\veps}$,
\begin{eqnarray*}
\|\wtil\alpha_j(\wtil x_{k_\veps})\kern1pt\what T^{n_j}\wtil y
-\what x\|_{\what\X}
&\kern-6pt=\kern-6pt&
\|(\wtil\alpha_j(\wtil x_{k_\veps})\kern1pt
(\what T|_{\wtil\X})^{n_j}\wtil y-\wtil x_{k_\veps})
+(\wtil x_{k_\veps}-\what x)\|_{\what\X}                                 \\
&\kern-6pt\le\kern-6pt&
\|\wtil\alpha_j(\wtil x_{k_\veps})\kern1pt\wtil T^{n_j}\wtil y
-\wtil x_{k_\veps}\|_{\what\X}
+\|\wtil x_{k_\veps}\!-\what x\|_{\what\X}
\le\veps
\end{eqnarray*}
whenever ${j\ge j_{\kern1pt\veps,k_\veps}}.$ Take an arbitrary integer
${i\ge1}$ and set ${\veps=\frac{1}{i}}.$ Consequently, set
${k(i)=k_\veps=k_{\frac{1}{i}}}$ and
$j(i)=j_{\kern1pt\veps,k_\veps}\!=j_{\kern1pt\frac{1}{i}\kern-.5pt,k{(i)}}\!$
for every ${i\ge1}.$ Thus for every integer ${i\ge1}$ there is an integer
${j(i)\ge0}$ such that, by setting
$\wtil\alpha_j(\wtil x_{k(i)})=\wtil\alpha_j(\wtil x_{k_\veps})
=\wtil\alpha_j(\wtil x_{k_\veps}(\what x))$,
we get
${\|\wtil\alpha_j(\wtil x_{k(i)})\kern1pt\what T^{n_j}\wtil y
-\what x\|_{\what\X}}\le\smallfrac{1}{i}$
whenever ${j\ge j(i)}.$ Therfore, by taking ${j=j(i)}$,
$$
\|\wtil\alpha_{j(i)}(\wtil x_{k(i)})\kern1pt\what T^{n_{j(i)}}\wtil y
-\what x\|_{\what\X}\le\smallfrac{1}{i}
\quad\;\hbox{for every integer ${i\ge1}$}.
$$
So for the sequence $\{\wtil\alpha_{j(i)}(\wtil x_{k(i)})\}_{i\ge1}$ we get
${\|\wtil\alpha_{j(i)}
(\wtil x_{k(i)})\kern1pt\what T^{n_{j(i)}}\wtil y-\what x\|_{\what\X}\to0}$.
Hence
$$
\wtil\alpha_{j(i)}
(\wtil x_{k(i)})\kern1pt\what T^{n_{j(i)}}\wtil y\conv\what x.
$$
For each integer ${i\ge1}$, set
${\what\alpha_i(\what x)=\wtil\alpha_{j(i)}(\wtil x_{k(i)})
=\wtil\alpha_{j(i)}(\wtil x_{k(i)}(\what x))}$
and ${\what T^{n_i}\kern-1pt=\what T^{n_{j(i)}}}.$ Thus for every
${\what x\in\what\X}$ there exist a sequence of nonzero numbers
$\{\what\alpha_i(\what x)\}_{i\ge1}$ and a sequence
$\{\what T^{n_i}\}_{i\ge1}$ with entries from $\{\what T^n\}_{n\ge0}$ such
that ($**$) holds$.$ So ${\wtil y\in\wtil\X\sse\what\X}$ is a supercyclic
vector for ${\what T\!\in\B[\wtil\X]}.$ This proves (2), concluding the
proof \hbox{of Claim 1$.\!\!\!\qed$}

\vskip6pt\noi
But $T$ is an isometry on the normed space $({\X,\|\cdot\|_\vphi})$, and so
is its extension $\what T$ on the completion $({\what\X,\|\cdot\|_{\what\X}})$
of $({\X,\|\cdot\|_\vphi}).$ Thus, by (a) and (b), the Banach-space
isometry $\what T$ has a supercyclic vector$.$ This, however, leads to a
contradiction, since
$$
\hbox{\it an isometry on a Banach space is never supercyclic}\/
$$
\cite[Proof of Theorem 2.1]{AB} (see also \cite[Lemma 4.1(b)]{KD1}$\kern1pt).$
Outcome$:$ there is no super\-cyclic vector ${y\in\X}$ for the operator $T$
acting on the normed space $({\X,\|\cdot\|}\kern.5pt)$. \qed

\vskip6pt\noi
{\bf Theorem 5.2.}$\,$\cite{AB}
{\it If a power bounded operator\/ $T$ on a normed space\/ $\X$ is
supercyclic, then it is strongly stable}\/.

\vskip6pt\noi
{\it Proof}\/$.$
Suppose $T$ is a power bounded operator and set $\beta=\sup_n\|T^n\|>0$.

\vskip6pt\noi
{\it Claim II\/}$.$
If ${y\in\X}$ is supercyclic for $T$, then ${T^ny\to0}$.

\vskip6pt\noi
{\it Proof of Claim II}\/$.$
Let ${0\ne y\in\X}$ be supercyclic for $T\kern-1pt.$ Equivalently, let
${y\in\X}$ be a supercyclic unit vector (${\|y\|=1}$) for $T\kern-1pt.$
Suppose ${T^ny\not\to0}.$ Thus by Claim 0 (cf$.$ proof of Theorem 5.1)
$0<\inf_n\|T^ny\|$; that is, there exists a $\delta>0$ for which
$$
\delta<\|T^ny\|                                                \eqno{\rm(i)}
$$
for every integer $n.$ Moreover, since $T$ is power bounded and supercyclic,
it follows from Theorem 5.1 that $T$ is not of class $C_{1{\textstyle\cdot}}.$
So there exists a unit vector ${v\in\X}$ (${\|v\|=1}$) such that ${T^nv\to0}.$
But since $y$ is a supercyclic vector for $T\kern-1pt$, there exists a
sequence $\{\alpha_k\}_{k\ge0}$ of nonzero numbers and a sequence
$\{T^{n_k}\}_{k\ge0}$ with entries from $\{T^n\}_{n\ge0}$ such that
${\alpha_kT^{n_k}y\to v}$, and hence there exists a positive integer $k_0$
such that if ${k\ge k_0}$, then $\|{\alpha_kT^{n_k}y-v}\|<\frac{1}{2}.$ So
$|\,{\|\alpha_kT^{n_k}y\|-1}\,|=|\,{\|\alpha_kT^{n_k}y\|-\|v\|}\,|
\le{\|\alpha_kT^{n_k}y-v\|}<\frac{1}{2}$,
which implies ${1-\|\alpha_kT^{n_k}y\|}<\frac{1}{2}$, and therefore
$\frac{1}{2}=1-\frac{1}{2}<$ $\|\alpha_kT^{n_k}y\|
\le|\alpha_k|\sup_n\|T^n\|\|y\|=\beta|\alpha_k|.$
Hence for ${k\ge k_0}$
$$
\smallfrac{1}{2\beta}<|\alpha_k|.                             \eqno{\rm(ii)}
$$
Now take any $\gamma$ for which
${0<\kern-.5pt\gamma\kern-1pt<\kern-.5pt\frac{\delta}{2\beta^2}}.$ Since
${\alpha_kT^{n_k}y\to v}$, there is ${k_1\kern-1pt\ge k_0}$ such that if
${k\ge k_1}$, then $\|\alpha_kT^{n_k}y-v\|<\smallfrac{\gamma}{2}.$ Since
$\|{\alpha_kT^{n_k+m}y\kern-.5pt-\kern-1ptT^mv}\|
\le\|T^m\|\kern.5pt{\|\alpha_kT^{n_k}y\kern-.5pt-v\|}$,
we get for ${k\ge k_1\ge k_0}$ and ${m\ge0}$,
$$
\|\alpha_kT^{n_k+m}y-T^mv\|<\beta\smallfrac{\gamma}{2}.      \eqno{\rm(iii)}
$$
Next, as ${T^nv\to0}$, take $m$ sufficiently large such that
$$
\|T^mv\|<\beta\smallfrac{\gamma}{2}.                          \eqno{\rm(vi)}
$$
Thus according to the above four displayed inequalities, for $k$ and $m$
large enough
\begin{eqnarray*}
\smallfrac{\delta}{2\beta}
&\kern-6pt<\kern-6pt&
\smallfrac{\|T^{n_k+m}y\|}{2\beta}
<\|T^{n_k+m}y\|\,|\alpha_k|=\|\alpha_kT^{n_k+m}y\|                         \\
&\kern-6pt\le\kern-6pt&
\|\alpha_kT^{n_k+m}y-T^mv\|+\|T^mv\|
<\beta\smallfrac{\gamma}{2}+\beta\smallfrac{\gamma}{2}
=\beta\gamma<\smallfrac{\delta}{2\beta},
\end{eqnarray*}
which is a contradiction$.$ Hence if ${y\in\X}$ is supercyclic for
$T\kern-1pt$, then ${T^ny\to0}$, which concludes the proof of Claim II. $\qed$

\vskip6pt\noi
Suppose $T$ is supercyclic$.$ As is readily verified, every nonzero vector in
the projective orbit $\Oe_T([y])$ of any supercyclic vector $y$ for $T$ is
again a supercyclic vector for $T$ (see, e.g., \cite[Lemma 5.1]{Kub2})$.$
Hence
$$
\hbox{\it the set of all supercyclic vectors for $T$ is dense in the normed
space $\X.$}
$$
Therefore, according to Claim II, ${T^ny\to0}$ for every $y$ in a dense subset
of $\X.$ Thus, since $T$ is power bounded, this implies that ${T^nx\to0}$ for
every ${x\in\X}.$ Outcome$:$ $T$ is strongly stable.                    \qed

\vskip6pt
Such detailed proofs of Theorems 5.1 and 5.2 are fundamental to supporting the
forthcoming discussion in Section 6.
\vskip-1pt\noi

\section*{{\bf 6}. Weak Supercyclicity and Weak Stability}

$\!$The problem considered in this section focuses on only three forms of
supercyclici\-ty, viz., (strong) supercyclicity, weak l-sequential
supercyclicity, and weak supercy\-clicity$.$ So consider the following
implications taken from the diagram in \hbox{Section 3}:
\vskip5pt\noi
\centerline{
{\esc supercyclic}
$\;\limply\;$
{\esc weakly l-sequentially supercyclic}
$\;\limply\;$
{\esc weakly supercyclic},
}
\vskip5pt\noi
where the converses fail$.$ In fact, it was exhibited in
\cite[Example 3.6]{BM1} a weakly supercyclic unitary operator on a Hilbert
space, which actually is weakly l-sequentially supercyclic (cf$.$
\cite[p.10, Corollary to Example 3.6]{BM1}), as we saw in Subsection 4.3$.$
Thus, since there is no unitary supercyclic operator, because there is no
supercyclic isometry on a \hbox{Banach} space, as we saw in the proof of
Theorem 5.1 (cf$.$ \cite[Proof of Theorem 2,1]{AB} or
\cite[Lemma 4.1(b)]{KD1}$\kern.5pt$), it follows that weak l-sequential
supercyclicity does not imply supercyclicity$.$ Also, as we saw in
Subsection 4.5, it was shown in \cite[Corollary 1.3]{Shk1} that weak
supercyclicity does not imply weak l-sequential supercyclicity --- by
exhibiting a weakly supercyclic unitary operator that is not weakly
sequentially super\-cyclic, and so it is not weakly l-sequentially
supercyclic.

\vskip6pt
Since supercyclicity implies strong stability, as we saw in Theorem 5.2, the
next question arises quite naturally$.$ For a power bounded operator on a
normed space,
$$
\hbox{\it does weak supercyclicity imply weak stability}\,?
$$
This is a nontrivial question that remains unanswered$.$ Pacing P\'olia
\cite[p.114]{Pol} when he says ``can you imagine a more accessible related
problem'', we might replace weak supercyclicity in the above question by weak
l-sequential supercyclicity$.$ This again leads to another nontrivial question
that also seems to remain unanswered so far$.$ Regarding such a new question
inquiring whether weak l-sequential supercyclicity implies weak stability, the
following diagram summarises the above line of reasoning.
\vskip4pt\noi
$$
\matrix{
&\hbox{\esc supercyclicity} &\limply& \hbox{\esc strong stability}
\phantom{\Big|}                                                          \cr
& \big\Downarrow            &       & \big\Downarrow
                                                                         \cr
&\hbox{\esc weak l-sequential supercyclicity}
                            &\query & \hbox{\esc weak stability}.
\phantom{\Big|}                                                          \cr}
$$

\vskip4pt
There are, however, some affirmative answers for particular classes of
operators$.$ For instance, when weak l-sequential supercyclicity coincides
with supercyclicity for a given class of operators, as is the case, for
example, of compact operators$.$ In these cases, Theorem 5.2 ensures a
positive answer to the question for such classes$.$ The connection between
weak l-sequential supercyclicity and weak stability for the case when $\X$ is
a Radon--Riesz space was explored in \cite{KD1}, \cite{KV1}, and \cite{KD5}$.$
Recall that a Radon--Riesz space is a normed space $\X$ for which an
$\X$-valued sequence $\{x_n\}$ converges strongly if and only if it converges
weakly and the sequence of norms $\{\|x_n\|\}$ converges to the norm of the
limit (see, e.g., \cite[Definition 2.5.26]{Meg})$.$ Hilbert spaces are
Radon--Riesz spaces (cf$.$ \cite[Problem 20, p.13]{Hal}).
 
\vskip6pt
Consider the proofs of Theorems 5.1 and 5.2$.$ In an attempt to bring the
arguments therein, concerning norm convergence, to weak convergence as far as
possible, the following results have been established in
\cite[Theorem 3.1]{KD4}.

\vskip6pt\noi
{\bf Theorem 6.1.}$\,$\cite{KD4}
{\it Let\/ $T\!$ be a power bounded operator of class\/
$C_{1{\textstyle\cdot}\!}$ on a normed space\/ ${(\X,\|\cdot\|)}.$ Then
\vskip2pt\noi
\begin{description}
\item{$\kern-8pt$\rm(a)$\kern2pt$}
there is a norm\/ ${\|\cdot\|_\vphi}$ on\/ $\X$ for which\/ $T$ is an
isometry on\/ ${(\X,\|\cdot\|_\vphi)}$.
\vskip4pt\noi
\item{$\kern-8pt$\rm(b)$\kern2pt$}
If a vector\/ ${y\in\X}$ is weakly l-sequentially supercyclic for\/ $T$ when
it acts on the normed space\/ ${(\X,\|\cdot\|)}$, then the same vector\/
${y\in\X}$ is weakly l-sequentially supercyclic for\/ $T$ when it acts on the
normed space\/ ${(\X,\|\cdot\|_\vphi)}$.
\vskip4pt\noi
\item{$\kern-8pt$\rm(c)$\kern2pt$}
If\/ $T$ is weakly l-sequentially supercyclic when acting on\/
${(\X,\|\cdot\|)}$, then it has an extension\/ $\what T$ on the
completion\/ ${(\what\X,\|\cdot\|_{\what\X})}$ of\/ ${(\X,\|\cdot\|_\vphi)}$
which is a weakly l-sequentially supercyclic isometric isomorphism\/.
\vskip4pt\noi
\item{$\kern-8pt$\rm(d)$\kern2pt$}
If\/ $T$ on\/ ${(\X,\|\cdot\|)}$ is weakly stable, then\/ $\what T$ on\/
${(\what\X,\|\cdot\|_{\what\X})}$ is weakly stable\/.
\end{description}
}

\vskip5pt
By comparing the statements in Theorem 6.1(a,b,c) with the proof of
Theorem 5.1, we can gauge the extent to which the proof of Theorem 5.1 can
survive the move from (strong) supercyclicity to weak l-sequential
supercyclicity.

\vskip5pt
A crucial argument in the proof of Theorem 5.1 is the fact that {\it there
is no super\-cyclic isometry on a Banach space}\/$.\!$ This was necessary
for concluding that {\it super\-cyclicity implies strong stability}\/$.$
Differently from the argument in the proof of The\-orem 5.1, there is
no counterpart to such a fact in the context of Theorem 6.1$.$ Indeed, as we
have seen in Subsection 4.3, {\it there are weakly l-sequentially supercyclic
unitary operators on a Hilbert space}\/$.$ A comparison between the proofs of
Theorem 6.1 as given in \cite[Theorem 3.1]{KD4} and that of Theorem 5.1 shows
that both proofs share the same structure, although the weak counterpart of
Theorem 5.2 (i.e., weak sta\-bility being implied by weak l-sequential
supercyclicity) is not \hbox{reached this way.}

\vskip5pt
A way to smooth down weak stability leads to the notion of weak
quasistability$.$ Recall that an operator $T$ on a normed space $\X$ is weakly
stable if $\lim_n\!|f(T^nx)|\kern-1pt=\kern-1pt0$

\goodbreak\noi
for every ${x\in\X}$, for
every ${f\in\X^*}\!.$ It is {\it weakly quasi\-stable}\/ if,
$$
{\liminf}_n|f(T^nx)|=0
\;\;\hbox{for every}\;
x\in\X,
\;\hbox{for every}\;
f\in\X^*\!.
$$
It is clear that weak stability implies weak quasistability, but the converse
fails \cite[Proposition 4.3]{KV2}$.$ Although the previous question remains
open, the next result (which is from \cite[Corollary 4.3]{KD5}$\kern.5pt$)
says that, for power bounded operators,
\vskip6pt\noi
\centerline{
{\esc weak l-sequential supercyclicity}
$\;\limply\;$
{\esc weak quasistability}.
}
\vskip1pt\noi

\vskip6pt\noi
{\bf Theorem 6.2.}$\,$\cite{KD5}
{\it Every weakly l-sequentially supercyclic power bounded operator on a
normed space is weakly quasistable}\/.

\vskip6pt
This shows that such a mild version of weak stability (viz., weak
quasistability) is implied by weak l-sequential supercyclicity, as strong
stability was implied by (strong) supercyclicity in Theorem 5.2$.\kern-1pt$
So weak l-sequential supercyclicity ensures that the origin is a weak limit
point of the orbit of every vector \cite[\hbox{Corollary 4.4}]{KD5} --- and
so the origin is a weak accumulation point of it (cf$.$ Subsection 3.1):

\vskip6pt\noi
{\bf Corollary 6.3.}$\,$\cite{KD5}
{\it If\/ $T$ is a power bounded operator on a normed space\/ $\X$, then}\/
$$
\Oe_T([y])^{-wl}=\X
\;\;\hbox{\it for some}\;\;
y\in\X
\quad\limply\quad
0\in\Oe_T(x)^{-wl}
\;\;\hbox{\it for every}\;\;
x\in\X.
$$
\vskip-4pt\noi

\vskip6pt
Let ${y\in\X}$ be a weakly l-sequentially supercyclic vector for
${T\kern-1pt\in\BX}.\kern-1pt$ Take an arbitrary
${x\in\X\backslash\Oe([y])}.\kern-1pt$ Then there exist a scalar sequence
$\{\alpha_k\}$ and a subsequence $\{T^{n_k}\}$ of $\{T^n\}$ for which
${\alpha_kT^{n_k}y\wconv x}.$ These are referred to as a scalar sequence
$\{\alpha_k\}$ and a subsequence $\{T^{n_k}\}$ {\it of weak l-sequential
supercyclicity}\/ (of\/ $T\kern-1pt$ for $x$ with respect to $y.)$ If
$\{T^{n_k}\}$ is such that ${\sup_k(n_{k+1}-n_k)<\infty}$, then it is called
{\it a boundedly spaced}\/ subsequence of $\{T^n\}.\kern-1pt$ Theorem 6.4
below is from \cite[Theorems 6.2 and 6.3]{KV2}, where item (a) is a
consequence of Theorem 6.2~above.

\vskip6pt\noi
{\bf Theorem 6.4.}\cite{KV2}
{\it Let\/ $\kern-1ptT$ be weakly l-sequentially supercyclic operator on a
normed space}\/~$\X$.
\begin{description}
\item{$\kern-6pt$\rm(a)$\kern2pt$}
{\it If\/ $T\kern-2pt$ is power bounded, and if there is a boundedly spaced
subsequence\/ $\{T^{n_k}\}$ of weak l-sequential supercyclicity of\/ $T$ for
each vector\/ ${x\in\X\backslash\Oe_T([y])}$, for \hbox{every} weakly
l-sequentially supercyclic vector\/ ${y\kern-.5pt\in\kern-.5pt\X}$, then\/
$T\kern-1pt$ is weakly \hbox{stable}}\/.
\vskip4pt\noi
\item{$\kern-6pt$\rm(b)$\kern1pt$}
{\it If\/ $T$ is weakly stable, then all scalar sequences\/ $\{\alpha_k\}$ of
weak l-sequential su\-per\-cyclicity of\/ $T\kern-1pt$ are unbounded$.$ In
other words, $\{\alpha_k\}$ is unbounded for \hbox{every}
${x\in\X\backslash\Oe_T([y])}$ and every weakly l-sequential supercyclic
vector}\/ ${y\in\X}$.
\end{description}
\vskip-2pt

\vskip6pt
Theorems 6.1 and 6.2 lead to descriptions of weakly l-sequentially supercyclic
isometries on a Hilbert space$.\kern-1pt$ The result below is from
\cite[Corollary 4.3]{KD4} as an application of a inner-product-space version
of Theorem 6.1 given in \cite[\hbox{Theorem 3.2}]{KD4}.

\vskip6pt\noi
{\bf Corollary 6.5.}$\,$\cite{KD4}
{\it On a Hilbert space, if a weakly l-sequentially supercyclic operator is
similar to an isometry, then it is similar to an unitary operator}\/.

\vskip6pt
This ensures that {\it if a weakly l-sequentially supercyclic operator on a
Hilbert space is power bounded and power bounded below, then it is similar to
unitary operator}\/ \cite[Remark 41]{KD4} --- we will return to this topic
in Theorem 7.1(c) in the next section$.$ On the other hand, as an application
of Theorem 6.2 and of the upcoming Theorem 6.7, the next result (from
\cite[Corollary 5.5]{KD5}$\kern.5pt$) sets a suitable starting point, although
not in a chronological order, for our discussion towards a characterisation of
weakly l-sequentially supercyclic unitary operators.
\goodbreak\noi

\vskip6pt\noi
{\bf Corollary 6.6.}$\,$\cite{KD5}
{\it If an isometry on a Hilbert space is weakly l-sequentially super\-cyclic,
then it is a weakly quasistable singular-continuous unitary operator}\/.

\vskip3pt
Being weakly quasistable, we may inquire, in general,
$$
\hbox{\it when is a unitary operator weakly stable}\,?
$$

\vskip0pt
Recall that (i) a scalar spectral measure for a normal operator is a positive
finite measure equivalent to its spectral measure; (ii) a normal operator on a
separable Hilbert space has a scalar spectral measure; (iii) any form of
cyclicity implies separability --- thus assume separability; (iv) a unitary
operator is absolutely continuous, singular-discrete, or singular-continuous
if its scalar spectral measure is absolutely continuous, singular-discrete, or
singular-continuous with respect to normalised Lebesgue measure on the
$\sigma$-algebra of Borel subsets of the unit circle, respectively; and
(v) every unitary operator is the orthogonal direct sum of an absolutely
continuous unitary, a singular-discrete unitary, and a singular-continuous
unitary (where any part may be missing)$.$ Also recall that (vi) an absolutely
continuous unitary is always weakly stable; (vii) a singular-discrete unitary
is never weakly stable; and (viii) there exist singular-continuous unitaries
that are either weakly stable or weakly unstable (see, e.g.,
\cite[Section 3]{Kub1})$.\kern-1pt$ Therefore,
\vskip4pt\noi
{\narrower\narrower
{\it a unitary operator is weakly stable if and only if its
singular-\hbox{continuous part}\/$\kern-11pt$
$($if it exists\/$)$ is weakly stable and its singular-discrete part does not
exist}\/.
\vskip0pt}

\vskip6pt
Examples of weakly stable and weakly unstable singular-continuous unitary
operators were given in \cite[Proposition 3.2 and Proposition 3.3]{Kub1}$.$
Also, along the lines discussed in \cite[pp.10,12]{BM1} (see \cite[p.1383]{BM}
as well) and \cite[Proposition 4.1]{Kub1}, it was shown in
\cite[Theorem 4.2]{Kub1} that singular-continuous unitary operators partially
characterise weakly l-sequentially supercyclic unitary operators.

\vskip6pt\noi
{\bf Theorem 6.7.}$\,$\cite{Kub1}
{\it Every weakly l-sequentially supercyclic unitary operator is
sin\-gular-continuous}\/.

\vskip6pt
Theorem 6.7 prompts the question inquiring about the converse:
\vskip4pt\noi
{\narrower\narrower
{\it is every singular-continuous unitary operator on a separable Hilbert
space weakly l-sequentially supercyclic}\,?
\vskip0pt}

\vskip6pt\noi
Maybe not, but if yes, then we would get a power bounded weakly
l-sequentially supercyclic (unitary) operator that is not weakly stable (but
weakly quasistable~by Theorem 6.2) since, as we saw above, there exist
singular-continuous unitary opera\-tors that are not weakly stable$.$ The next
question also remains unanswered$.$ Actually, in spite of Theorem 6.2, or
because of it, and having in mind Theorem 5.2, the question below that closes
this section is nearly the same one that has opened~it.

\vskip6pt\noi
{\bf Question 6.8.}$\,$\cite{KD1}
{\it Is a weakly l-sequentially supercyclic\/ power bounded operator weakly
stabile}\,?

\vskip6pt
This was posed in \cite[p.61]{KD1}$.$ Similarly (see, e.g.,
\cite[Question 5.2]{Kub1}),
$$
\hbox{\it is every weakly\/ l-sequentially supercyclic unitary operator
weakly stable}\,?
$$

\section*{{\bf 7}. Spectrum of Weakly {\rm l}-Sequentially Supercyclic
Operators}

We close the paper by giving a description of the spectrum of a weakly
l-se\-quentially supercyclic operator.

\vskip6pt
Take the spectrum $\sigma(T)$ of an operator $T\kern-1pt$ acting on a complex
Banach \hbox{space $\X\kern-1pt$}, which is nonempty and compact, and consider
its classical partition, namely,
$\sigma(T)=
{\sigma_{\kern-1ptP}(T)\cup\sigma_{\kern-1ptR}(T)\cup\sigma_{\kern-.5ptC}(T)}$,
where $\sigma_{\kern-1ptP}(T)$ is the point spectrum, $\sigma_{\kern-1ptR}(T)$
is the residual spectrum, and $\sigma_{\kern-.5ptC}(T)$ is the continuous
spectrum$.$ Let $\DD$ and $\TT$ be the open unit disk and the unit circle in
the complex \hbox{plane $\CC.$} For any subset $\Lambda$ of the complex plane
$\CC$, set $\Lambda^*\!=\{\lambda\in\CC\!:\overline\lambda\in\Lambda\}$, the
set of all complex conjugates of all \hbox{elements from $\Lambda$}.

\vskip6pt
Perhaps the spectral saga of hypercyclicity has its origins in 1982, when
Kitai \cite[Theorems 2.7 and 2.8]{Kit} has shown that if a Banach-space
operator $T$ is hyper\-cyclic, then ${\sigma(T)\cap\TT\ne\void}$, as well as
every component of $\sigma(T)$ meets $\TT\kern-.5pt.\kern-.5pt$ This has
triggered the next results on weak hypercyclicity and, consequently, on weak
supercyclicity$.$ As briefly mentioned in Section 4 (Subsections 4.1 and 4.2),
concerning weak hyper\-cyclicity, Pr\v ajitur\v a \cite[Theorem 2.2]{Pra} has
shown in 2005 that $\sigma_{\kern-1ptP}(T^*)=\void$ for a weakly hypercyclic
operator $T$ on a Hilbert space$.$ This extended a previous result of Herero
\cite[Proposition 2.2]{Her} in 1991, where such an expression had been
established for hypercyclic operators --- see also
\cite[Corollary 2.4]{Kit}$.$ It was shown by Dilworth and Troitsky
\cite[Theorem 3]{DT} in 2003 that the expression ${\sigma(T)\cap\TT\ne\void}$
survives for a weakly hypercyclic operator $T$ on a Banach space.

\vskip6pt
The above paragraph dealt with spectral properties of hypercyclic and weakly
hypercyclic operators$.$ Regarding supercyclicity, it was verified
that ${\#\sigma_{\kern-1ptP}(T^*)\le1}$ (i.e., the cardinality of
$\sigma_{\kern-1ptP}(T^*)$ is not greater than one) for supercyclic operators
on a Hilbert space by Herrero \cite[Proposition 3.1]{Her} in 1991, which was
extended to operators on a normed space by Ansari and Bourdon
\cite[Theorem 3.2]{AB}) in 1997, and further extended to supercyclic operators
on a locally convex space by Peris \cite[Lemma 1, Theorem 4]{Per} in 2001$.$
But the weak topology on a normed space is a locally convex subtopology of the
locally convex norm topology (see, e.g.,
\cite[Theorems 2.5.2 and 2.2.3]{Meg})$.$ Then the latter extension holds in
particular on a normed space under the weak topology, thus including weakly
supercyclic operators, and so weakly l-sequentially supercyclic operators, on
normed spaces$.$ This has also been obtained for weakly supercyclic operators
on a Hilbert space by Pr\v ajitur\v a \cite[Theorem 3.2]{Pra} in 2005$.$ Under
additional assumptions of power boundedness and being of class
$C_{1{\textstyle\cdot}}$, it was shown by Kubrusly and Duggal
\cite[Corollary 4.1]{KD4} in 2021 (as a corollary of Theorem 6.1) that
${\#\sigma_{\kern-1ptP}(T^*)=\#\sigma_{\kern-1ptP}(T)=0}$ for every weakly
l-sequentially super\-cyclic operators acting on Hilbert spaces.

\vskip6pt
We will proceed along these lines in Theorem 7.1 below for weakly
l-sequentially supercyclic operators on a Hilbert space, but first we need
the following result.

\vskip6pt
{\it An invertible power bounded operator on a Hilbert space, whose inverse
also is power bounded, is similar to a unitary operator}\/$.$ This is a
classical result due to \hbox{Sz$.$ Nagy} \cite[Theorem 1]{Nag} that will be
required in the proof of the next theorem, which is similarly stated as
follows$:$ {\it a Hilbert-space operator is similar to a unitary operator if
and only if it is power bounded and power bounded below}\/$.$ A normed space
version for operators similar to an isometry was later given by Koehler and
Rosenthal \cite[Theorem 2]{KR}$:$ {\it a normed-space operator is similar to
an isometry if and only if it is power bounded and power bounded below}\/.

\vskip6pt
$\!$Recall from Theorem 5.1$:$ a power bounded operator of class
$C_{1{\textstyle\cdot}}\kern-.5pt$ is not \hbox{supercyclic}.

\vskip6pt\noi
{\bf Theorem 7.1.}$\kern-.5pt$
{\it Let\/ $\kern-.5ptT\kern-.5pt$ be a power bounded operator of class\/
$C_{1{\textstyle\cdot}}\kern-1pt$ on a \hbox{Hilbert space\/
$\kern-.5pt\H\kern-.5pt.$}}
\begin{description}
\item{$\kern-4pt$(a)}
{\it If\/ $T$ is weakly l-sequentially supercyclic, then}
$$
\sigma(T)=\sigma_{\kern-.5ptC}(T)\sse\DD^-\!
\quad\hbox{\it and}\quad
\sigma(T)\cap\TT\ne\void.
$$
\item{$\kern-4pt$(b)}
{\it If\/ $T$ is weakly l-sequentially supercyclic with closed range, then}
$$
\sigma(T)=\sigma_{\kern-.5ptC}(T)\sse\DD^-\backslash B_\veps(0)
\quad\hbox{\it and}\quad
\sigma(T)\cap\TT\ne\void.
$$
{\it for some open ball\/ $B_\veps(0)$ with radius\/ ${\veps\in(0,1]}$
centred at the origin of}\/ $\H.$
\vskip6pt\noi
\item{$\kern-4pt$(c)}
{\it If\/ $T$ is weakly l-sequentially supercyclic and power bounded below,
then\/ $T$ is similar to a unitary operator and}
$$
\sigma(T)=\sigma_{\kern-.5ptC}(T)\sse\TT.
$$
\item{$\kern-4pt$(d)}
{\it If\/ $T$ is weakly l-sequentially supercyclic isometry, then\/ $T$ is a
singular continuous unitary operator\/, so that
$$
\sigma(T)=\sigma_{\kern-.5ptC}(T)\sse\TT
$$
is the support of a finite positive measure that is singular-continuous with
respect to normalised Lebesgue measure on the\/ $\sigma$-algebra of Borel
\hbox{subsets of\/ $\TT$\kern-.5pt.}}
\end{description}

\vskip6pt\noi
{\it Proof}\/$.$
(a)
If $T$ is power bounded on a complex Banach space, then
${\sigma(T)\sse\DD^-\!}$ (cf$.$ Section 2)$.$ It was shown in
\cite[Corollary 4.1]{KD4} (by using an inner-product-space version of
Theorem 6.1) that if a power bounded operator $T\kern-1pt$ of class
$C_{1{\textstyle\cdot}}\kern-.5pt$ on a com\-plex Hilbert space is weakly
l-sequentially supercyclic, then
$$
\sigma_{\kern-1ptP}(T)=\sigma_{\kern-1ptP}(T^*)=\void,
$$
where $T^*$ stands for the Hilbert-space adjoint of $T.$ Since for
Hilbert-space operators
$$
\sigma_{\kern-1ptR}(T)
=\sigma_{\kern-1ptP}(T^*)^*\backslash\sigma_{\kern-1ptP}(T),
$$
it follows from the above two displayed expressions that
$$
\sigma_{\kern-1ptP}(T)=\sigma_{\kern-1ptR}(T)=\void,
$$
and so $\sigma(T)=\sigma_{\kern-.5ptC}(T)\sse\DD^-\!.$ Moreover, as $T$ is not
strongly stable (because it is of class $C_{1{\textstyle\cdot}}$), it is not
uniformly stable, and so its spectral radius is not strictly less than one
(see Section 2 again), which implies that its spectrum intersects the unit
circle (i.e., ${\sigma(T)\cap\TT\ne\void}$, since
${\sigma(T)\kern-1pt\sse\kern-1pt\DD^-}).$ This concludes the proof of
\hbox{item (a)}.

\vskip6pt\noi
(b)
The range of a weakly l-sequentially supercyclic operator is dense (see, e.g.,
\cite[Lem\-ma 5.3]{KD5}$\kern,5pt).$ Thus $T$ is surjective if it has a closed
range$.$ Since $T$ is injective (as ${0\not\in\sigma_{\kern-1ptP}(T)}$), it
is invertible with a bounded inverse, and so ${0\not\in\sigma(T)}$, which
implies ${\sigma(T)\sse\DD^-\backslash B_\veps(0)}$ since $\sigma(T)$ is a
closed set included in the closed unit \hbox{disk by (a).}

\vskip6pt\noi
(c)
If, besides being power bounded, $T$ is power bounded below, then $T$ is
similar to a unitary operator (as we saw above)$.$ As similarity preserves the
spectrum, and the spectrum of a unitary operator is included in the unit
circle, we get \hbox{(c) from (a)}.

\vskip6pt\noi
(d)
By Corollary 6.6, every weakly l-sequentially supercyclic isometry on a
Hilbert space is a singular-continuous unitary$.$ Recall that an isometry is
trivially power bounded of class $C_{1{\textstyle\cdot}}$, and a unitary
operator is, in addition, power bounded below$.$ As the spectrum of a
singular-continuous unitary operator is the support of its singular-continuous
spectral measure, the result in (c) comes from (b). \qed

\section*{Acknowledgment}

I thank an anonymous referee who made sensible observations throughout the
text, contributing to the improvement of the paper.

\bibliographystyle{amsplain}

\end{document}